\documentclass[preprint,12pt]{elsarticle}
\usepackage{latexsym,amssymb,amsmath,fancyhdr}
\usepackage{todonotes}
\usepackage{color}
\usepackage{palatino}
\usepackage{a4wide}
\usepackage{graphicx}
\usepackage{amsfonts} 
\usepackage{amsthm}
\usepackage[TABBOTCAP]{subfigure}
\usepackage[latin1]{inputenc}
\usepackage[bookmarks=true,colorlinks=true,backref]{hyperref}

\newcommand{\nOne}[0]{^{n\scalebox{0.8}{\!+\!}1}} 
\newcommand{\kOne}[0]{_{k\scalebox{0.8}{\!+\!}1}} 
\newcommand{\nHalf}[0]{^{n\scalebox{0.8}{\!+\!}\frac{1}{2}}} 
\newcommand{\mygamma}[0]{M}

\begin{document}

\begin{frontmatter}
\title{Time integration for diffuse interface models for two-phase flow}


\author{Sebastian Aland}
\address{Institut f\"ur wissenschaftliches Rechnen, TU Dresden, 01062 Dresden, Germany}
\ead{sebastian.aland@tu-dresden.de}

\begin{abstract}
We propose a variant of the $\theta$-scheme for diffuse interface models for two-phase flow, 
together with three new linearization techniques for the surface tension.
These involve either additional stabilizing force terms, or a fully implicit coupling of the Navier-Stokes and Cahn-Hilliard equation.

In the common case that the equations for interface and flow are coupled explicitly, we find a time step restriction 
which is very different to other two-phase flow models and in particular is independent of the grid size.
We also show that the proposed stabilization techniques can lift this time step restriction.

Even more pronounced is the performance of the proposed fully implicit scheme which is stable for arbitrarily large time steps.
We demonstrate in a Taylor flow application that this superior coupling between flow and interface equation can render diffuse interface models even computationally cheaper and faster than sharp interface models.

\end{abstract}

\begin{keyword}
time integration, diffuse interface model, dominant surface tension, time stability, CFL condition, Navier Stokes, Cahn Hilliard, linearization
\end{keyword}

\end{frontmatter}

\tableofcontents

\section{Introduction}\label{sec:introduction}

The numerical simulation of two-phase flows has reached some importance in microfluidic applications.
In the last decade, diffuse interface (or phase-field) models have become a valuable alternative to the more established sharp interface methods (e.g. Level-Set, Arbitrary Lagrangian-Eulerian, Volume-Of-Fluid). 
The advantages of diffuse interface methods include the possibility to easily handle moving
contact lines and topological transitions as well as the fact that they do not require any reinitialization or convection stabilization. 
The corresponding equations involve a Navier-Stokes (NS) equation coupled to a convective Cahn-Hilliard (CH) equation.
A lot of efficient spacial discretization techniques and solvers for these equations have been proposed (e.g. \cite{KayWelford_SIAMJSC_2007}).
However, not much work has been done on time integration strategies and efficient coupling between the NS and the CH equation, which we will address in this paper. 

But at first, let us introduce the diffuse interface method more carefully. 
The method was originally developed to model solid-liquid phase transitions, see e.g. \cite{Andersonetal_ARFM_1998,Emmerich_AP_2008,SingerSinger_RPP_2008}. The interface thereby is represented as a thin layer of finite thickness and an auxiliary function, the so-called phase field, is used to indicate the phases. The phase field function varies smoothly between distinct values in both phases and the interface can be associated with an intermediate level set of the phase field function. Diffuse interface approaches for mixtures of two immiscible, incompressible fluids lead to the NS-CH equations and have been considered by several authors, see e.g. \cite{Jacqmin_JCP_1999,Feng_SIAMJNA_2006,KayWelford_SIAMJSC_2007,Dingetal_JCP_2007}. The simplest model reads: 
\begin{eqnarray}
\rho(c) \left( \partial_t {\bf u} + ( {\bf u} \cdot \nabla) {\bf u} \right) &=& 
-\nabla {p} + \nabla \cdot \left(\nu(c) \mathbf{D}(\mathbf{u})\right)+ {\bf F} +  \mu \nabla c , 
\label{nsch1}\\
\nabla \cdot {\bf u} &=& 0, 
\label{diff incompr}\\
\partial_t c + {\bf u} \cdot \nabla c &=& \nabla\cdot\left( M(c) \nabla \mu \right), \\
\mu &=& {\tilde\sigma}\epsilon^{-1} W'(c)-{\tilde\sigma}\epsilon\Delta c, \label{nsch4}
\end{eqnarray}
in the domain $\Omega$. Here ${\bf u}$, ${p}$, $c$ and $\mu$ are the velocity, pressure, phase field variable and chemical potential, respectively. 
The function $W(c)$ is a double well potential, here we use $W = 1/4 (c^2-1)^2$ which ensures that $c\approx \pm 1$ in two fluid phases, respectively. 

The function $M(c)$ is a mobility, $\epsilon$ defines a length scale over which the interface is smeared out. 
In general for the diffuse interface fluid method, it is desirable to keep $M$ small such that one primarily gets advection.
At the same time the mobility needs to be big enough to ensure that the interface profile stays accurately modeled
and the interface thickness is approximately constant. 
Furthermore ${\bf D}({\bf u})=\nabla{\bf u} + \nabla{\bf u}^T$ is the strain tensor, $\rho(c)$, $\nu(c)$ and $\mathbf{F}$ are the (phase dependent) density, viscosity and body force. The parameter $\tilde\sigma$ is a scaled surface tension which is related to the physical surface tension by $\tilde{\sigma} = \sigma \frac{3}{2\sqrt{2}}$.			
There are efficient solvers available to discretize and solve the Eqs. (\ref{nsch1})-(\ref{nsch4}) in space (see e.g. \cite{KayWelford_SIAMJSC_2007}).

Surface tension is a major component of all multiphase fluid models and hence
various spatial discretizations of the surface tension force for diffuse-interface models have been proposed
(e.g. \cite{Kim_JCP_2005}).
The surface tension force $\mu\nabla c$ introduces a strong coupling between the NS equation providing the flow field and the CH equation evolving the phase field. 
This is very similar to sharp interface models for two-phase flow where the same interface-to-flow coupling introduces a severe time step restriction
 of the form \cite{Badalassi_JCP_2003, Brackbilletal_JCP_1992}:
 \begin{align}
   \tau < C \rho^{1/2} h^{3/2}\sigma^{-1/2}.
 \end{align}
Here, $\tau$ is the maximum time step size, $\rho$ the average density of both fluids and $h$ the grid size.
The above CFL-like restriction is particularly strong for large effective surface tensions, e.g. when small physical length scales are considered.
It is usually assumed that this restriction also holds for diffuse interface models (e.g. in \cite{KimLowengrub_IFB_2005}). In Sec. \ref{sec:stability} we will show that this assumption is wrong.

\begin{figure}[h]
\centering
\includegraphics[width=0.8\textwidth]{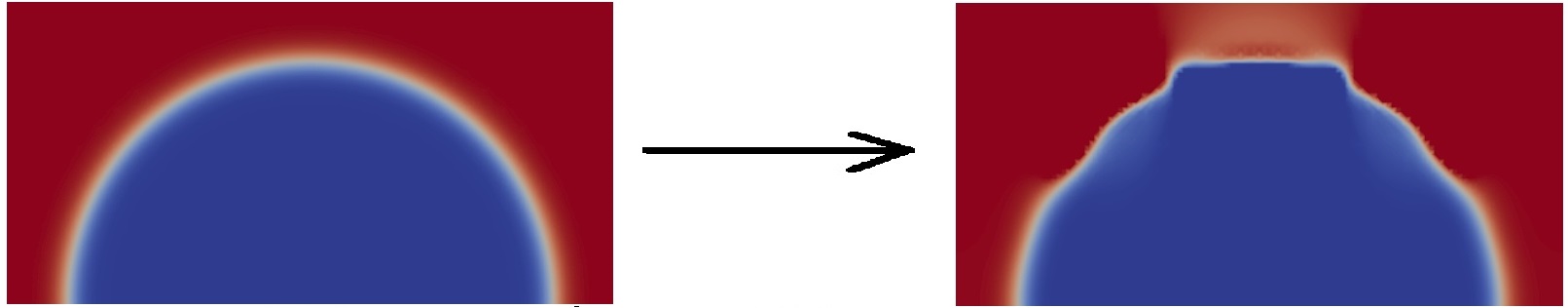}
\caption{Evolution of a semi-circular bubble under the diffuse interface model with too large time steps.
Left: initial shape; Right: Evolved shape after ten time steps.
}
\label{fig:surface_tension}
\end{figure}

However, also for diffuse interface models there is some time step restriction which can make computations extremely costly, even in cases when the interface is supposed to hardly move.
Fig. \ref{fig:surface_tension} shows such a case of a perfectly circular interface, which is almost stationary. 
However, if too big time steps  are chosen,  even such an equilibrated surface will start to wobble and finally break up.
In sharp interface models, there are techniques to overcome such time step restrictions \cite{Hysing_IJNMF_2006}. 
To the best of the authors' knowledge there is no such technique available for diffuse interface models yet.
We will develop techniques to improve the coupling between the NS and the CH equations, which will turn out to lift the time step restrictions significantly.

Apart from increasing the computational performance, there is a second reason to develop better time integration schemes for diffuse interface models.
The simple time discretization schemes available often imply the need to stabilize the system by choosing a relatively high CH mobility.
But this high artificial diffusion perturbs the simulation results since matched asymptotic analysis shows the convergence of diffuse-interface methods toward the sharp interface equations only for small CH mobility \cite{Abelsetal_WS_2012}.
Therefore better time integration strategies would not only speed-up the simulations but also
allow to take smaller (more physical) CH mobility and thus improve the accuracy of diffuse-interface methods. 

The structure of the remaining paper is as follows. 
Secs. \ref{sec:time discretization} and \ref{sec:linearization and coupling} will introduce a simple variant of the $\theta$ scheme as well as a block Gauss-Seidel coupling strategy. The main attention is given to Sec. \ref{sec:advanced linearization} where some new improved coupling techniques for diffuse interface models are presented. The solution of the resulting systems is discussed in Sec. \ref{sec:space discretization}.
In Sec. \ref{sec:results} we perform numerical tests. In particular a CFL-condition for diffuse interface methods is numerically derived and it is shown that the new proposed coupling methods can, for some problems, result in an extreme gain of performance. 
Finally, conclusions are drawn in Sec. \ref{sec:conclusions}.

\section{Time discretization: A variant of the $\theta$-scheme}\label{sec:time discretization}

In this section we adopt the well-known $\theta$-scheme for the time discretization of the NS-CH equations. 
Let the time interval $[0,T]$ be divided in $N$ subintervals of size $\tau^n$, $n=1...N$.
We define the discrete time derivative of a (solution) variable $v$ to be $d_t v\nOne:=(v\nOne-v^{n})/\tau^n$, where the upper index denotes the time step number.
For a shorter notation we introduce 
\begin{align}
g({\bf u}, c, \mu)&:= -\rho(c) {\bf u} \cdot \nabla {\bf u} +\nabla \cdot \left(\nu(c) \mathbf{D}(\mathbf{u})\right)+ {\bf F}(c) + \mu\nabla c,\\
f({\bf u}, c, \mu)&:= -{\bf u} \cdot \nabla c + \nabla\cdot(M(c)\nabla\mu).
\end{align}
For a constant $\theta\in[0,1]$ we propose the following variant of the $\theta$-scheme:
\begin{align}
\rho\nHalf  \partial_t {\bf u}\nOne + \nabla p\nOne=& 
\theta g({\bf u}^{n{+}1}, c\nOne, \mu\nOne) + (1-\theta)g({\bf u}^{n}, c^{n}, \mu^{n}), 
\label{nsch_theta1}\\
\nabla \cdot {\bf u}\nOne =& 0, 
\label{nsch_theta2}\\
\partial_t c\nOne =& \theta f({\bf u}\nOne, c\nOne, \mu\nOne) + (1-\theta)f({\bf u}^{n}, c^{n}, \mu^{n}) , \label{nsch_theta3}\\
\mu\nOne =& {\tilde\sigma}\epsilon^{-1} W'(c\nOne)-{\tilde\sigma}\epsilon\Delta c\nOne, \label{nsch_theta4}
\end{align}
where $\rho\nHalf = (\rho(c\nOne)+\rho(c^n))/2$ denotes an approximation to the  intermediate densitiy.
Note that the system (\ref{nsch_theta1})-(\ref{nsch_theta4}) differs from the standard $\theta$-scheme by the use of $\rho\nHalf$. To derive a standard $\theta$-scheme one would have to divide the NS equation by $\rho(c)$. Then the $\theta$-scheme can be applied to the system $\partial_t {\bf u} = [g({\bf u}, c,\mu)-\nabla p]/\rho(c)$. This leads to more complicated equations requiring much more implementation effort. The above variant circumvents this problem by  using $\rho\nHalf$. One easily verifies that the method has the same consistency order and stability properties as the standard $\theta$-scheme.
To be more precise, the method is of second order if $\theta=0.5$ and of first order if $\theta\neq 0.5$. The most important cases are $\theta=1$ (backward Euler) and $\theta=0.5$ (Crank-Nicolson). Both are A-stable and are often used to discretize two-phase flow problems in practice. 
We will also restrict our numerical experiments to these two cases. The biggest disadvantage of the Crank-Nicolson scheme is that it has no smoothing property, i.e. it does not smooth high frequencies. Therefore it is in some cases appropriate to use the optimally smoothing backward Euler scheme, although it has lower order.

\section{Linearization and coupling}\label{sec:linearization and coupling}
The first question that arises when looking at equations (\ref{nsch_theta1})-(\ref{nsch_theta4}) is: How to couple the NS and CH equations.
The nonlinear coupling between ${\bf u}\nOne, c\nOne$ and $\mu\nOne$ can be treated by several decoupling strategies. 
One usually applies an iterative strategy where each time step requires multiple solves of the governing equations until the approximation error in the solution variables is sufficiently small.
To our knowledge, only sequential coupling of the NS and CH equations has been considered in the literature. Hence, both equations are solved separately, using the solution of the other equation explicitly from a previous computation.

\subsection{Block Gauss-Seidel decoupling} \label{sec:block gauss seidel}
The simplest decoupling strategy is the block Gauss-Seidel method.
This method is widely used, e.g. in \cite{Dingetal_JCP_2007, GruenKlingbeil_arXiv_2012}.
We use subscript indices to denote the variables of the sub-iteration while superscripts still denote the time step number. Then the block Gauss-Seidel strategy applied in every time step reads:
\begin{enumerate}
\item Initialize the sub-iteration with the values from the last time step:
\begin{align}
{\bf u}_0 = {\bf u}^n,~~~~ c_0=c^n,~~~~ \mu_0 = \mu^n
\end{align}
\item for k=0,1,...\\
$\bullet$ Solve the NS equation to get ${\bf u}\kOne$ and $p\kOne$:
\begin{align}
\rho\nHalf  \partial_t {\bf u}\kOne + \nabla p\kOne=& 
\theta g({\bf u}\kOne, c_k, \mu_k) + (1-\theta)g({\bf u}^{n}, c^{n}, \mu^{n}), 
\label{nsch_GS}\\
\nabla \cdot {\bf u}\kOne =& 0,\label{nsch_GS2}
\end{align}
where $\rho\nHalf$ is approximated by $(\rho(c_k)+\rho(c^n))/2$. \\
$\bullet$ Solve the CH equation to get $c\kOne$ and $\mu\kOne$:
\begin{align} 
\partial_t c\kOne =& \theta f({\bf u}\kOne, c\kOne, \mu\kOne) + (1-\theta)f({\bf u}^{n}, c^{n}, \mu^{n}) , \label{nsch_GS3}\\
\mu\kOne =& {\tilde\sigma}\epsilon^{-1} W'(c\kOne)-{\tilde\sigma}\epsilon\Delta c\kOne, \label{nsch_GS4}
\end{align}
$\bullet$ proceed to the next $k$
\item stop the iterative process when a given tolerance is reached, e.g. $||c_{k+1}-c_k||<tol$, and set the variables at the new time step to
\begin{align}
{\bf u}\nOne = {\bf u}\kOne,~~~~ c\nOne=c\kOne,~~~~ \mu\nOne  = \mu\kOne
\end{align}
\end{enumerate}
Note that in Eq. (\ref{nsch_GS}) only previously calculated instances of the CH variables $c$ and $\mu$ appear.
One can interpret the Gauss-Seidel strategy as a fix point iteration where the fix point operator consists of one solve of the NS equation and one solve of the CH equation. The fix point iteration consists in applying this operator on the iterates $({\bf u}_k, c_k, \mu_k)$ multiple times until convergence is reached. Note, that convergence can be very slow due to the high nonlinearity of the operators. In fact a contraction is only assured for very small time steps and divergence may occur if the time step is choosen too large.

\subsection{Linearization}\label{sec:simple linearization}
There are three remaining nonlinear terms: ${\bf u}\kOne\cdot\nabla{\bf u}\kOne$ in Eq. (\ref{nsch_GS}), $W'(c\kOne)$ in Eq. (\ref{nsch_GS4}) and $\nabla\cdot(M(c\kOne)\nabla\mu\kOne)$ in Eq. (\ref{nsch_GS3}). 
Using a Taylor series expansion gives the second order linear approximations: 
\begin{align}
W'(c\kOne)&\approx W'(c_k)+W''(c_k)(c\kOne-c_k) \label{linear_W}\\
{\bf u}\kOne\cdot\nabla{\bf u}\kOne &\approx {\bf u}_k\cdot\nabla{\bf u}\kOne + {\bf u}\kOne\cdot\nabla{\bf u}_k - {\bf u}_k\cdot\nabla{\bf u}_k \label{linear_inertia}\\
\nabla\cdot(M(c\kOne)\nabla\mu\kOne) &\approx \nabla\cdot(M(c_k)\nabla\mu\kOne)+ \nabla\cdot(M'(c_k)(c\kOne-c_k)\nabla\mu_k) \label{linear_mobility}
\end{align}
If small length scales are considered, as in most applications of diffuse interface models, convection does not dominate the NS equation.
Hence, a simpler linearization than Eq. (\ref{linear_inertia}) is sufficient. Here we skip the last two terms on the RHS of Eq. (\ref{linear_inertia}).
Note that this does not affect the accuracy of the method and would at most slow down the convergence if the problem was dominated by convection.
Similarly, the last term of Eq. \eqref{linear_mobility} can be omitted since it has very little effect on the convergence speed, at least in all of our applications.

\subsection{Special case: semi-implicit Euler} \label{sec:semi implicit Euler}
A special case of the linearized $\theta$-scheme 
is when $\theta=1$ and only one iteration of Eqs. \eqref{nsch_GS}-\eqref{nsch_GS4} is performed.
The resulting method can be written as the following semi-implicit Euler scheme:
\begin{align}
\rho(c^n)  \partial_t {\bf u}\nOne + \nabla p\nOne=& 
 -\rho(c^n) {\bf u}^n \cdot \nabla {\bf u}^{n+1} +\nabla \cdot \left(\nu(c^n) \mathbf{D}(\mathbf{u}\nOne)\right)+ {\bf F}(c^n) + \mu^n\nabla c^n, 
\label{nsch_semi_implicit}\\
\nabla \cdot {\bf u}\nOne =& 0, 
\label{nsch_semi_implicit2}\\
\partial_t c\nOne =& -{\bf u}\nOne \cdot \nabla c\nOne + \nabla\cdot(M(c^n)\nabla\mu\nOne), \label{nsch_semi_implicit3}\\
\mu\nOne =& {\tilde\sigma}\epsilon^{-1} W'(c\nOne)-{\tilde\sigma}\epsilon\Delta c\nOne, \label{nsch_semi_implicit4}
\end{align}
where $W'(c\nOne)$ is again linearized by Eq. \eqref{linear_W}.
In each time step Eqs. \eqref{nsch_semi_implicit}-\eqref{nsch_semi_implicit2} and Eqs. \eqref{nsch_semi_implicit3}-\eqref{nsch_semi_implicit4} can be solve sequentially. 
This semi-implicit Euler scheme is the simplest time-stepping strategy. It is of first order accuracy but may still give good results if sufficiently small time steps are used.
The scheme is has been used by many authors \cite{KayWelford_SIAMJSC_2007, Aland2011_bijels} including a benchmark comparison of  diffuse interface with level-set and VOF methods which showed good agreement \cite{Aland_IJNMF_2011}.

\subsection{Defect correction scheme}
It is sometimes useful to apply a defect correction scheme when solving the fix-point iteration in Eqs. \eqref{nsch_GS}-\eqref{nsch_GS4}.
Therefore, the solution variable is split into its previous value plus an update value (here denoted with a star): $u\kOne = u_k + u^*$.
Now, the linear system is only solved for the update variable.
Hence, solving a system of the form $Au\kOne=b$ is replaced by solving $Au^*=b-Au^k$.
Left and right hand side of the latter equation are of the order of the defect $u^*$, which minimizes errors in computer arithmetic. 
Another advantage of this approach is that iterative solvers sometimes perform better, i.e. need less iterations, when looking for solutions of the defect corrected scheme.
To get a better approximation $u\kOne$, the update equation may be altered to ${\bf u}\kOne = {\bf u}_k + \omega{\bf u}^*$, where
the step length $\omega$ can be found by a line search strategy. 
The resulting Richardson type method is described in detail in \cite{Turek1999}.

Combining a defect correction scheme with the above block Gauss-Seidel splitting and the linearizations \eqref{linear_W}-\eqref{linear_mobility} is equal to applying one step of a Newton iteration alternately to each of the two subsystems (NS and CH).

\section{Advanced linearization techniques}\label{sec:advanced linearization}
We will now develop techniques to improve the coupling between the NS and the CH equations, which will turn out to lift the time step restrictions significantly.
In an iterative scheme, like the block Gauss-Seidel coupling, this is equivalent to linearizing the surface tension force more efficiently. 
We will start by coupling the NS and CH equation really implicitly by assembling both equations in one large system.

\subsection{Method 1: Fully-coupled scheme} \label{sec:implicit coupling}
The reason for the instability of the NSCH system for larger time steps or high surface tensions is the explicit coupling of the 
NS and CH equations. To be more precise it is the chemical potential that may oscillate (i.e. alters its sign) in 
every time step. Consequently, the best way to stabilize the system is by taking the chemical potential from the new time step, i.e. $\mu\kOne\nabla c_k$ instead of 
$\mu_k\nabla c_k$ in the NS equation, while still advecting the phase field with the new velocity ${\bf u}\kOne$. 
This can be done by assembling both, the NS and CH equations, in one large system. 
Hence, step 2 of the block Gauss-Seidel iteration is replaced by
\begin{itemize}
\item[2.] for k=0,1,...\\
$\bullet$ Solve the coupled NS-CH equation
\begin{align}
\rho\nHalf  \partial_t {\bf u}\kOne + \nabla p\kOne=& 
\theta g({\bf u}\kOne, c_k, \mu\kOne) + (1-\theta)g({\bf u}^{n}, c^{n}, \mu^{n}), 
\label{nsch_coupled}\\
\nabla \cdot {\bf u}\kOne =& 0,\label{nsch_coupled2} \\
\partial_t c\kOne =& \theta f({\bf u}\kOne, c_k, \mu\kOne) + (1-\theta)f({\bf u}^{n}, c^{n}, \mu^{n}) , \label{nsch_coupled3}\\
\mu\kOne =& {\tilde\sigma}\epsilon^{-1} W'(c\kOne)-{\tilde\sigma}\epsilon\Delta c\kOne, \label{nsch_coupled4}
\end{align}
$\bullet$ proceed to the next $k$
\end{itemize}
The slight difference to the block Gauss-Seidel iteration from Sec. \ref{sec:block gauss seidel} is the usage of of $c_k$ in Eq. \eqref{nsch_coupled3} and $\mu\kOne$ in Eq. \eqref{nsch_coupled}. 
We will see in Sec. \ref{sec:results} that this has a strong stabilizing effect, may speed up convergence of the fix-point scheme and allows larger time-steps.

The fix point iteration with this fully coupled system corresponds to a Newton iteration of the fully coupled system.
Hence one can expect this method in most cases to perform as well as a Newton iteration of the fully coupled equations, but without the need to compute the full Jacobian of the system.

We also tested other variants of assembling NS and CH in one system, for instance using the surface tension force $\mu_k\nabla c\kOne$.
We found all other variants not to give any improvements in stability. 
The key point here is really the use of the new curvature (contained in $\mu\kOne$). This is a big advantage of diffuse interface models over other two-phase flow models, because the curvature ($\mu$) here is a solution variable and can therefore be implicitly coupled to the NS equation.

\subsection{Method 2: Linearization of the chemical potential} \label{sec:S1}

In this section we present an alternative way which avoids solving NS and CH in one system. We will derive a stabilizing term which can easily be added to the NS equation.
From the previous section we know that it is desirable to replace the surface tension force $\theta\mu_k\nabla c_k $ occuring in Eq.(\ref{nsch_GS}) by the more implicit version
\begin{align}
F_{st}=\theta\mu\kOne\nabla c_k
\end{align}
Instead of using $\mu_k$ which is a first order approximation of the new chemical potential $\mu\kOne$, we will now derive a second order approximation.
The idea is to predict $\mu\kOne$ from the available variables $c_k$ and ${\bf u}\kOne$. To start, we use that the surface tension force $F_{st}$ can also be written as \cite{Feng_SIAMJNA_2006}: 
\begin{align}
F_{st}&=-\theta\tilde{\sigma}\epsilon\nabla\cdot(\nabla c\kOne\otimes\nabla c_k) \label{fst}.
\end{align}
Next, let us assume that the movement of $c$ is primarily driven by advection, i.e. the influence of the chemical potential in Eq.\eqref{nsch_GS3} is neglected. Substracting Eq.(\ref{nsch_GS3}) of step $k+1$ from the same equation in step $k$ gives
\begin{align}
c\kOne&=c_k-\tau\theta({\bf u}\kOne\nabla c\kOne-{\bf u}_k\nabla c_k)\\
&\approx c_k-\tau\theta({\bf u}\kOne-{\bf u}_k) \nabla c_k
\end{align}
Now, we may insert this expression in Eq.\eqref{fst} and get
\begin{align}
F_{st}&\approx -\theta\tilde{\sigma}\epsilon\nabla\cdot(\nabla c_k\otimes\nabla c_k) 
+\tau\theta^2\tilde{\sigma}\epsilon\nabla\cdot(\nabla(  
({\bf u}\kOne\!-\!{\bf u}_k) \nabla c_k
)\otimes\nabla c_k) \label{fst2}
\end{align}
The first term in Eq.\eqref{fst2} can be transformed back in $\theta\mu_k\nabla c_k$ which gives
\begin{align}
F_{st}\approx &~ \theta\mu_k\nabla c_k
+\tau\theta^2\tilde{\sigma}\epsilon\nabla\cdot(\nabla({\bf u}\kOne\!-\!{\bf u}_k)\cdot\nabla c_k\otimes\nabla c_k) \nonumber \\
&+\tau\theta^2\tilde{\sigma}\epsilon\nabla\cdot((\nabla\nabla c_k\cdot({\bf u}\kOne\!-\!{\bf u}_k))\otimes\nabla c_k) \label{fst3}
\end{align}
The last term in Eq.\eqref{fst3} turned out to be very small in all our simulations and did not influence neither accuracy nor convergence speed. We therefore omit it and get the stabilizing term 
\begin{align}
S_{1}&:=\tau\theta^2\tilde{\sigma}\epsilon\nabla\cdot(\nabla({\bf u}\kOne\!-\!{\bf u}_k)\cdot\nabla c_k\otimes\nabla c_k)\label{S_st}
\end{align}
which should be added to the RHS of the NS equation \eqref{nsch_GS}.
Note, that $S_1$ is a kind of Laplacian of the velocity field and can therefore be expected to have a stabilizing effect which will be confirmed in Sec. \ref{sec:results}.
Adding $S_{1}$ to Eq. \eqref{nsch_GS} does not affect the accuracy of the method since $S_1$ vanishes when the fix-point iteration converges. 
Also note, that the derivation assumed that the phase field $c$ is only advected, which corresponds to vanishing mobility. 
For larger mobilities $S_1$ should be scaled with an (unknown) factor $\omega\in[0,1]$. Here, we use $\omega=0.2$ which turned out to speed up the convergence of the fix-point method significantly.

\subsection{Method 3: Stabilizing surface Laplacian} \label{sec:S2}
The introduction of a stabilizing Laplacian of the velocity field in the previous section reminds of a very popular technique used in level-set methods first introduced by Dziuk \cite{Dziuk_NM_1990}. The idea is to use that $\kappa{\bf n} = -\Delta_\Gamma {id}_\Gamma$,  where $\kappa=\nabla\cdot{\bf n}$ is the mean curvature, ${\bf n}$ the normal and $\Delta_\Gamma \mathrm{id}$ the Laplace-Betrami of the identity mapping.
Hence, the implicit part of the surface tension force can be expressed as 
\begin{eqnarray}
F_{st} = - \sigma\theta \delta_k \Delta_{\Gamma_k} \mathrm{id}_{\Gamma_k}. \label{eq1}
\end{eqnarray}
with some surface Delta function $\delta_k$. 
To get a more implicit version of eq. (\ref{eq1}), one can replace $\mathrm{id}_{\Gamma_k}$ by $\mathrm{id}_{\Gamma\kOne}$. The latter can be approximated by
\begin{eqnarray}
 \mathrm{id}_{\Gamma\kOne} \approx ~~ x \rightarrow x+\theta\tau ({\bf u}\kOne-{\bf u}_k)(x). \label{id_approx}
\end{eqnarray}
which gives the surface tension force:
\begin{eqnarray} 
F_{st} &\approx& - \sigma\theta \delta_k \Delta_{\Gamma_k} \mathrm{id}_{\Gamma_k}  
- \sigma\theta^2 \delta_k\tau  \Delta_{\Gamma_k} ({\bf u}\kOne-{\bf u}_k) \\
&=&  - \sigma\theta \delta_k \Delta_{\Gamma_k} \mathrm{id}_{\Gamma_k}  
- \sigma\theta^2 \tau  \nabla_{\Gamma_k}\cdot\left(\delta_k \nabla_{\Gamma_k} ({\bf u}\kOne-{\bf u}_k)\right). \label{eq2}
\end{eqnarray}
The first term on the RHS corresponds to the surface tension force $\theta\mu_k\nabla c_k$ which is 	already included in our NS equation \eqref{nsch_GS}.
Consequently, the second term on the RHS of \eqref{eq2} is an additional stabilizing term which enters the NS equation. A diffuse interface version is given by
\begin{eqnarray} 
S_2:= - \sigma\theta^2\tau \nabla\cdot\left(|\nabla c_k| P  \nabla ({\bf u}\kOne-{\bf u}_k)\right) \label{S2}
\end{eqnarray}
where $P=\left(I-\frac{\nabla c_k\otimes\nabla c_k}{|\nabla c_k|^2}\right)$ is the surface projection.
Analogously to the derivation in the previous section, we assumed here (in Eq.\ref{id_approx}) that the interface is only advected.
Hence, the derivation only holds in the limit of vanishing mobility. For larger mobility we therefore scale $S_2$ in the same way as $S_1$ with a parameter $\omega\in[0,1]$.
Note, that adding $S_2$ does not affect the accuracy of the method, since it vanishes when the fix point iteration converges (and ${\bf u}\kOne\approx{\bf u}_k$).

\section{Space discretization and solvers}\label{sec:space discretization}

For the numerical solution of the partial differential equations we adapt existing algorithms for the NS-CH equation, e.g. \cite{VillanuevaAmberg_IJMP_2006,Do-QuangAmberg_PF_2009}. 
We use the finite element toolbox AMDiS \cite{amdis} for discretization with $P^2$ elements for ${\bf u}, c$ and $\mu$ and $P^1$ elements for the pressure $p$. 
In 2D we solve the resulting linear system of equations with UMFPACK \cite{umfpack}. 
For the larger systems in 3D we have to use preconditioned iterative solvers. Efficient preconditioners are available for the individual systems, that is when NS and CH are solved separately.
We apply the $F_p$ preconditioner \cite{KayWelford_SIAMJSC_2007} and an FGMRES iteration to solve the NS system.
For the CH system we use the preconditioner proposed in \cite{Boyanova_CMAM_2012} also with FGMRES iteration. 

It remains the case of the fully coupled NS-CH system, Eqs. \eqref{nsch_coupled}-\eqref{nsch_coupled4}, which has not been considered in the literature so far.
From the discretization of Eqs. \eqref{nsch_coupled}-\eqref{nsch_coupled4}  we obtain a system of the form
\begin{align}
\left[\begin{array}{cc}A_{NS}&M_c\\N_c&A_{CH}\end{array}\right] 
\left[\begin{array}{c} \vec{x} \\ \vec{y} \end{array}\right] 
=&
\left[\begin{array}{c} \vec{b_1} \\ \vec{b_2} \end{array}\right],  \label{system matrix}
\end{align}
where $\vec{x}$ contains the degrees of freedom of ${\bf u}$ and $p$, while $\vec{y}$ contains the degrees of freedom of $c$ and $\mu$.
$M_c$ and $ N_c$ denote the coupling terms between the NS and the CH system.
Our (simple) approach to solve this system is to combine the two preconditioners for the CH and NS systems.
Let $P_{NS}$ the $F_p$ preconditioner for NS and $P_{CH}$ the CH preconditioner. 
For the coupled system we use the matrix
\begin{align}
P = \left[\begin{array}{cc}P_{NS}&M_c\\&\frac{1}{\alpha} P_{CH}\end{array}\right]
\end{align}
as a right preconditioner, with a scaling factor $\alpha$. 
Hence, Eq. \eqref{system matrix} is replaced by solving the two systems
\begin{align}
\left[\begin{array}{cc}A_{NS}&M_c\\N_c& A_{CH}\end{array}\right]
\left[\begin{array}{cc}P_{NS}^{-1} &-\alpha P_{NS}^{-1}M_cP_{CH}^{-1}\\0&\alpha P_{CH}^{-1}\end{array}\right]
\left[\begin{array}{c} \vec{v}\\ \vec{w} \end{array}\right]
=&
\left[\begin{array}{c} \vec{b_1} \\ \vec{b_2} \end{array}\right]
\end{align}
and 
\begin{align}
\left[\begin{array}{c} \vec{x} \\ \vec{y} \end{array}\right]
=&
\left[\begin{array}{cc}P_{NS}^{-1} &-\alpha P_{NS}^{-1}M_cP_{CH}^{-1}\\0&\alpha P_{CH}^{-1}\end{array}\right]
\left[\begin{array}{c} \vec{v} \\ \vec{w} \end{array}\right].
\end{align}
A rigorous analysis of the matrix properties is still missing. However, in our numerical tests, this preconditioned system can be solved by an FGMRES iteration with $\alpha=0.1$.

\section{Numerical tests} \label{sec:results}
We now validate the proposed linearization schemes on different test scenarios.
First, we assess the numerical stability  of the time integration schemes (Sec. \ref{sec:stability}). 
Then we conduct a benchmark comparison of the different linearization methods (Sec. \ref{sec:benchmark}). 
At last we present an application to a Taylor-Flow simulation which clearly shows the superior performance 
of the proposed methods.
Throughout this section we use an iteration tolerance of $tol=10^{-10}$ (see Sec.\ref{sec:block gauss seidel}).
Furthermore, we terminate the block Gauss-Seidel scheme if no convergence is reached after 100 iterations. As 
stabilization constant we use $\omega=0.2$ (see Sec.\ref{sec:advanced linearization}).

\subsection{Stability investigations} \label{sec:stability}
In this section, we assess the time step stability of the proposed schemes.
Our goal is to find the maximum time step size at which a given two-phase flow configuration is solvable. 
We use configurations with different surface tensions $\sigma$, mobilities $\mygamma$, grid sizes $h$ and interface thicknesses $\epsilon$. 
In the end we want to find an estimate to predict the maximum time step size from these parameters.
As mentioned in Sec. \ref{sec:introduction} for other two-phase flow methods (e.g. Level-Set, ALE) the estimate gives the CFL-like condition
\begin{align}
\tau_{max} < C\rho^{1/2}\sigma^{-1/2} h^{3/2}
\label{cfl_standard} 
\end{align}
 with some non-dimensional constant $C$.
A common assumption is that such an estimate also holds for diffuse interface methods which we will contradict in the following. 

We use numerical testing to assess the time step stability, since analytical investigations on the coupled NS-CH system are extremely complicated. We restrict the numerical tests to the Crank-Nicolson scheme ($\theta=0.5$) \footnote{For $\theta=1$, the obtained maximum time step would be need to be divided by 2.}.
As initial condition we use ${\bf u} = 0$ and a phase field given by $c=\tanh((y-0.5)/\sqrt{2}\epsilon)$ in the domain $\Omega=[0,1]^2$ which represents a flat horizontal interface. Since the curvature is zero this initial condition corresponds to a stationary state.
To trigger an instability, a different random number $\in[-0.001,0.001]$ is added to each grid point of the phase field. 
Hence, the fix point iteration will not converge for too large time steps.
We then try to solve a single time step of the system with these initial conditions. 
We do this multiple times for varied time step sizes. We start with very small time steps which assure convergence.
As long as the system can be solved we increase the time step size (by a factor of 1.1) and start again.
At some point the system will not be solvable and we denote the corresponding time step as the maximum time step size $\tau_{\max}$ for this configuration.

\begin{table}
\centering
\begin{tabular}{l|rrrr}
			$\sigma$ & $10^2$ &$10^3$ &$10^4$ &$10^5$ \\\hline
			$\mygamma$ & $10^{-4}$& $10^{-5}$& $10^{-6}$& $10^{-7}$ \\\hline
			$\epsilon$ & $0.08$ & $0.04$ & $0.02$ & $0.01$ \\\hline
			$h$ & $2\epsilon$ & $\epsilon$ & $0.5 \epsilon$  \\\hline
 			$\rho$ & $10^{-1}$ &$10^0$ &$10^1$ &$10^2$ \\
\end{tabular}
\caption{Parameters used for stability estimations.}
\label{tab:stability_configurations}
\end{table}

We repeat this procedure for various numerical parameters shown in Tab. \ref{tab:stability_configurations}. We use equal densities in both phases and  choose $\nu=0.01$ to make sure that viscosity does not play a dominant role\footnote{analogously to the derivation of \eqref{cfl_standard}, see e.g.\cite{Brackbilletal_JCP_1992}}. 
Varying all other parameters ($\sigma,\rho,\epsilon,h,\mygamma$) independently gives a total number of $4\cdot 4\cdot 4\cdot 3\cdot 4=768$ test cases and
we denote the corresponding parameters for some test case $i\in\{1,...,768\}$ by subscripts: $\sigma_i,\rho_i,\epsilon_i,h_i,\mygamma_i$.
Our goal is to obtain a relationship between the test parameters and the corresponding maximally possible time step $\tau_{\max,i}$.
Assuming a multiplicative relationship with unknown exponents gives the nonlinear least-squares problem:
\begin{align}
\min_{\alpha} \sum_{i=1}^{768} \left| \log\left(\frac{\alpha_1 {h}_i^{\alpha_2} {\epsilon}_i^{\alpha_3} {\sigma}_i^{\alpha_4} {\mygamma}_i^{\alpha_5}\rho_i^{\alpha_6}}{\tau_{\max,i}} \right) \right|^2
\label{fitting}
\end{align}
where $\alpha$ contains the searched unknowns, $\alpha=(\alpha_1,...,\alpha_6)$. 
We use the routine \textsc{lscurvefit} in MatLab to solve \eqref{fitting} and obtain 
$\alpha = (7.603, 0.001, 0.916,  -0.326, 0.374, 0.673)$
 which is equivalent to the time step restriction 
\begin{align}
\tau_{\max} < 7.603 ~ h^{0.001} \epsilon^{ 0.916} \sigma^{-0.326} \mygamma^{0.374} \rho^{0.673}  \label{cfl_new1} 
\end{align}
The calculation of this time step restriction involved some rounding, in particular the maximum time steps may be over estimated by up to 10\%, since we increased them stepwise by 10\%.
These errors limit the precision of \eqref{cfl_new1} and justifies to round the obtained values. 
In this way, we get the following CFL-like time step restriction
\begin{align}
\tau_{\max} < 7.0~ \epsilon ~\sigma^{-1/3}~ \mygamma^{1/3}~ \rho^{2/3} \label{cfl_new2}
\end{align}
In strong constrast to \eqref{cfl_standard} the new time step restriction is independent of $h$. The reason for this lies in the fact that no sharp interface is used. The standard time step restriction \eqref{cfl_standard} is associated to the migration of capillary waves which might occur in sharp interface models with a wave length proportional to the grid size. In a diffuse interface context, the smallest wave length of capillary waves should be proportional to $\epsilon$ which would justify in \eqref{cfl_new2} to substitute one power of $h$ by $\epsilon$. Some additional smoothing is introduced by the CH diffusion which consequently also occurs in \eqref{cfl_new2}. 
Using that $\mygamma$ is measured in $m^3 s/kg$ we compute the physical unit of the RHS of \eqref{cfl_new2} and obtain s (seconds), which further justifies the new CFL-like condition.

\begin{figure}
\includegraphics[width=0.31\textwidth]{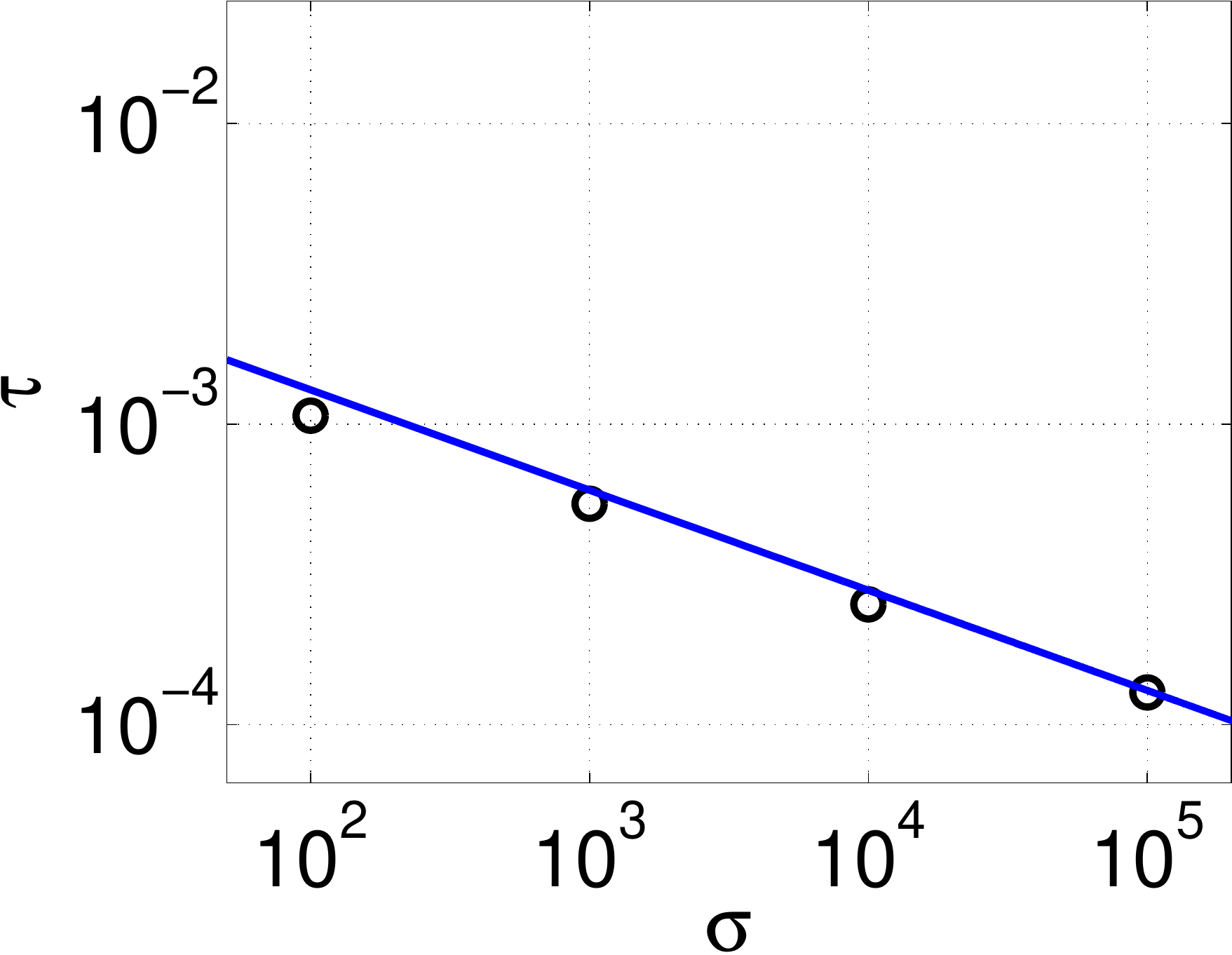}
\includegraphics[width=0.31\textwidth]{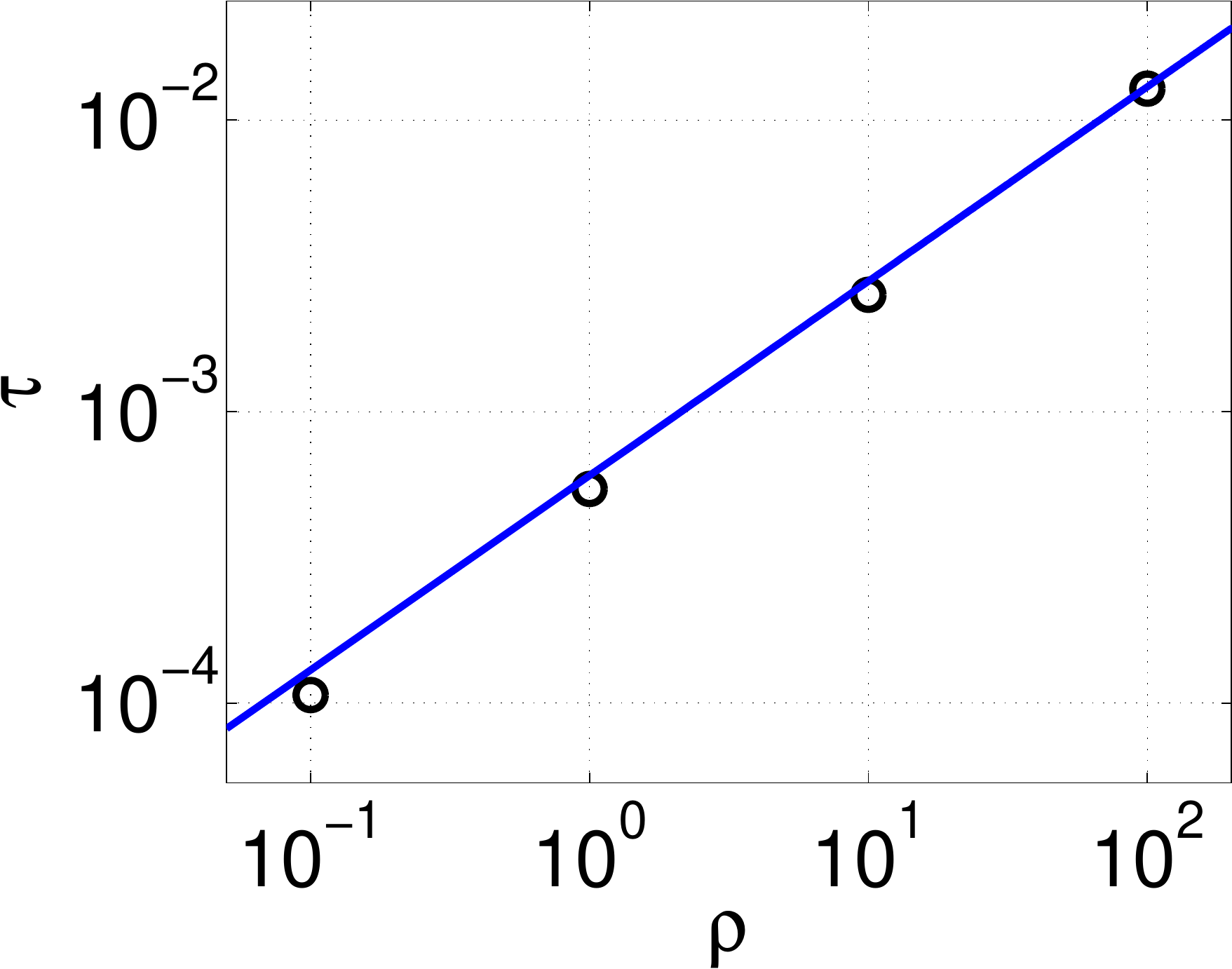}
\includegraphics[width=0.31\textwidth]{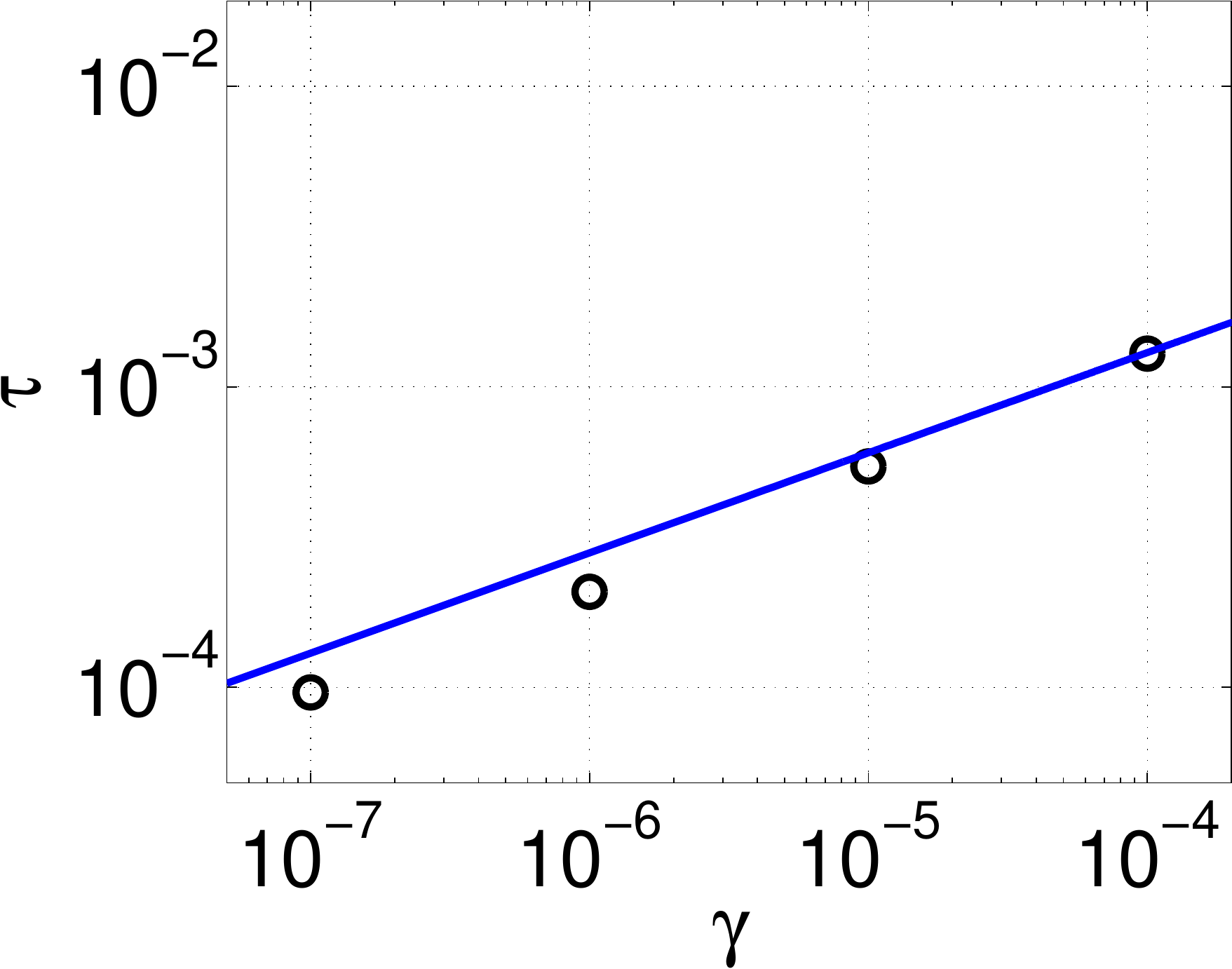}\\
\includegraphics[width=0.31\textwidth]{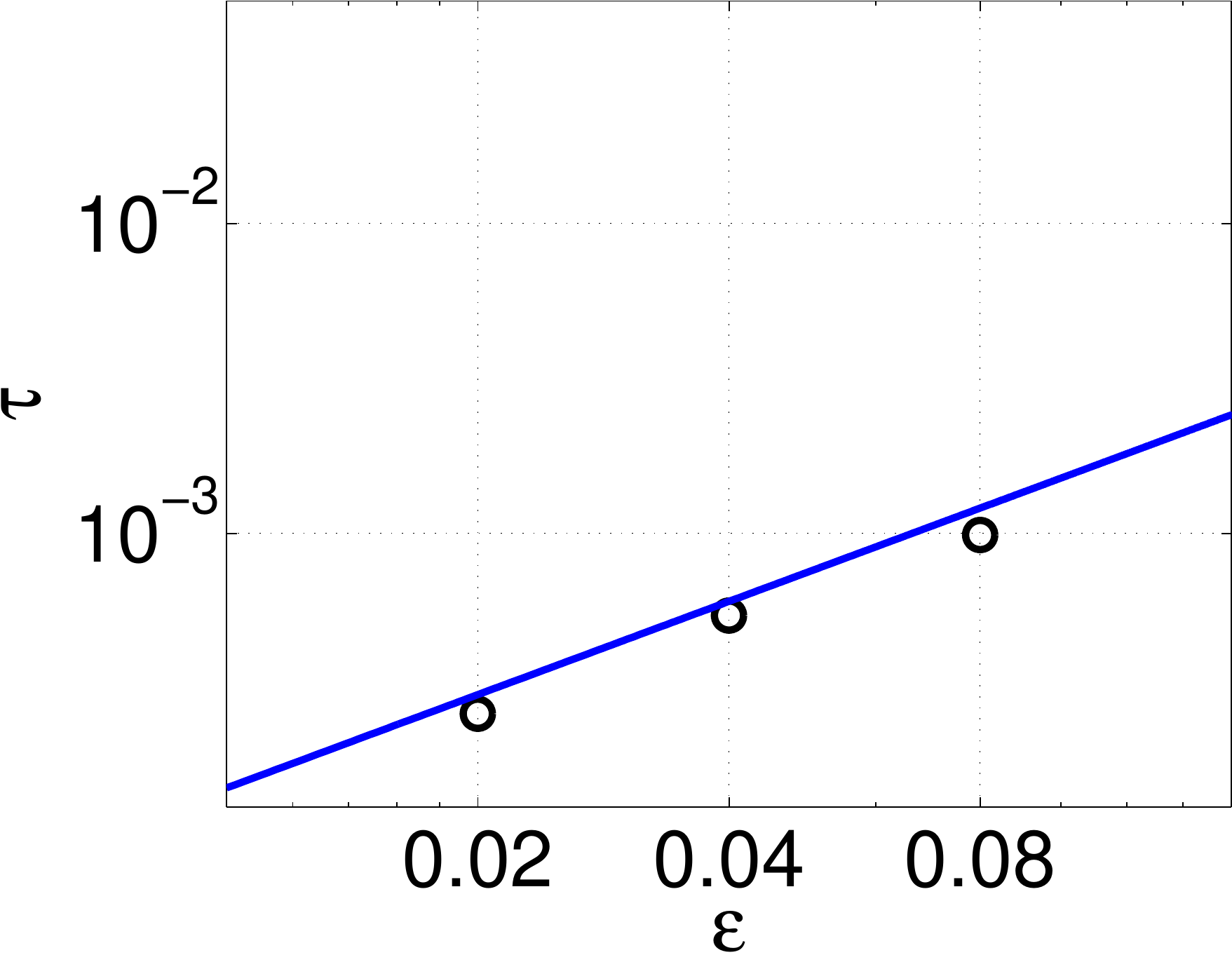}
\includegraphics[width=0.31\textwidth]{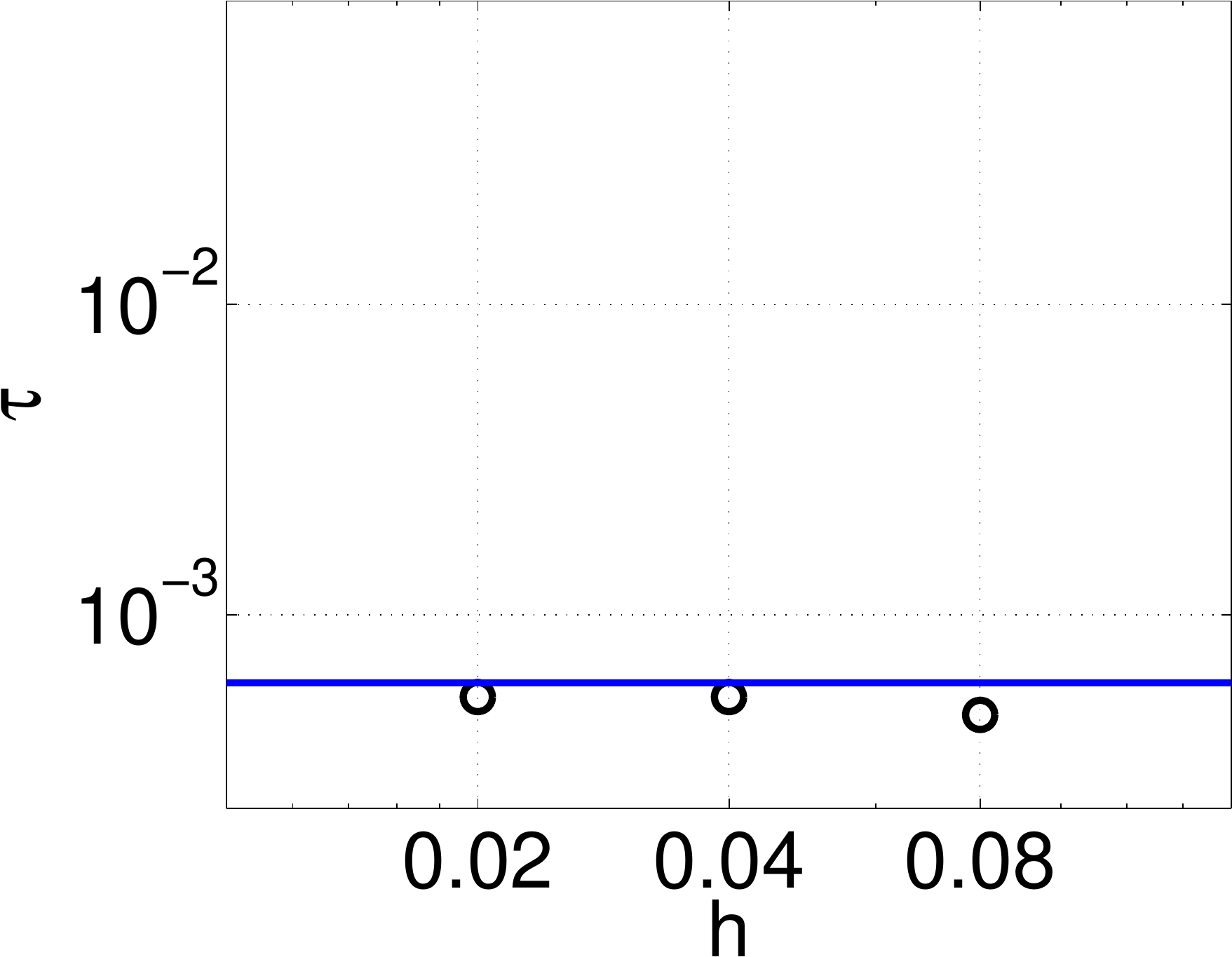}
\includegraphics[width=0.31\textwidth]{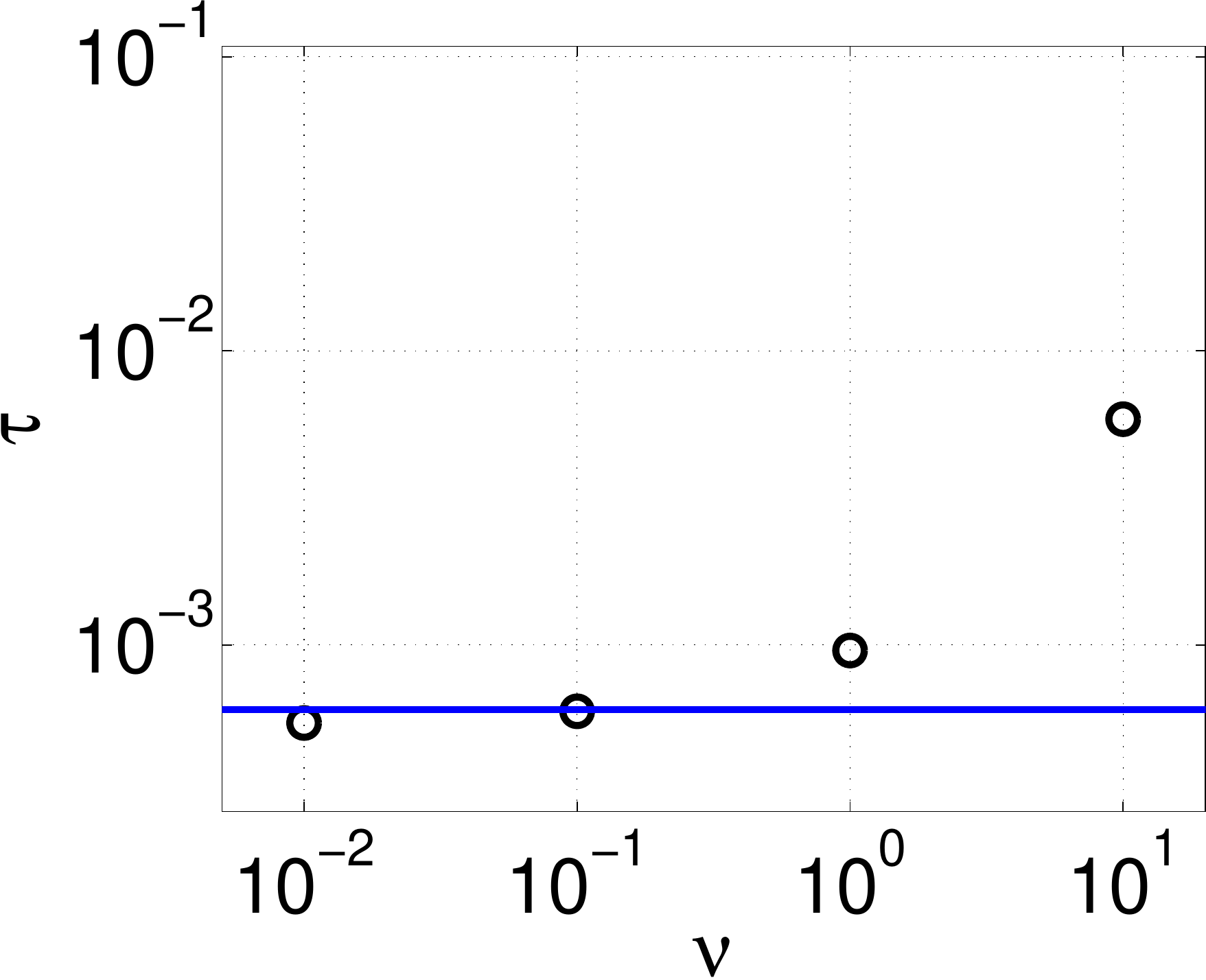}
\caption{Log-log plots comparing the maximum time step size with the CFL-condition \eqref{cfl_new2}. Either one of the variables $\rho,\sigma,\epsilon,\mygamma,h,\nu$ is varied while keeping the other variables fixed. }
\label{fig:stability}
\end{figure}

Figure \ref{fig:stability} shows the experimental maximum time step size compared to the CFL-condition \eqref{cfl_new2}.
Thereby, we vary one of the six variables ($\rho,\sigma,\epsilon,\mygamma,h,\nu$) while keeping the other variables fixed at $\rho=1.0,\sigma=10^3,\mygamma=10^{-5},\epsilon=0.04,h=0.04,\nu=0.01$. One can see an excellent agreement of the numerical data with the time step restriction curve. 
But also the limit of the derived CFL condition becomes apparent looking at the case of varied viscosity $\nu$ in Fig. \ref{fig:stability}.
According to our derivation, the time step restriction only holds in the limit of small viscosities (here $Re\lesssim 1$). 
Note, that this coincides with the famous CFL condition \eqref{cfl_standard} which was also derived for the small viscosity case \cite{Brackbilletal_JCP_1992}.  

We also assess the time step stability of the advanced linearization schemes proposed in Sec. \ref{sec:advanced linearization}. 
In particular the fully coupled sheme (Eqs. \eqref{nsch_coupled}-\eqref{nsch_coupled4}) shows a superior performance.
For all test cases the scheme converges in at most 10 iterations independently of the time step. Hence, for a configuration close to the stationary state the fully coupled scheme allows arbitrary large time steps. This outstanding property will further exploited in Sec. \ref{sec:taylor}. 

\begin{figure}
\centering
\includegraphics[width=0.4\textwidth]{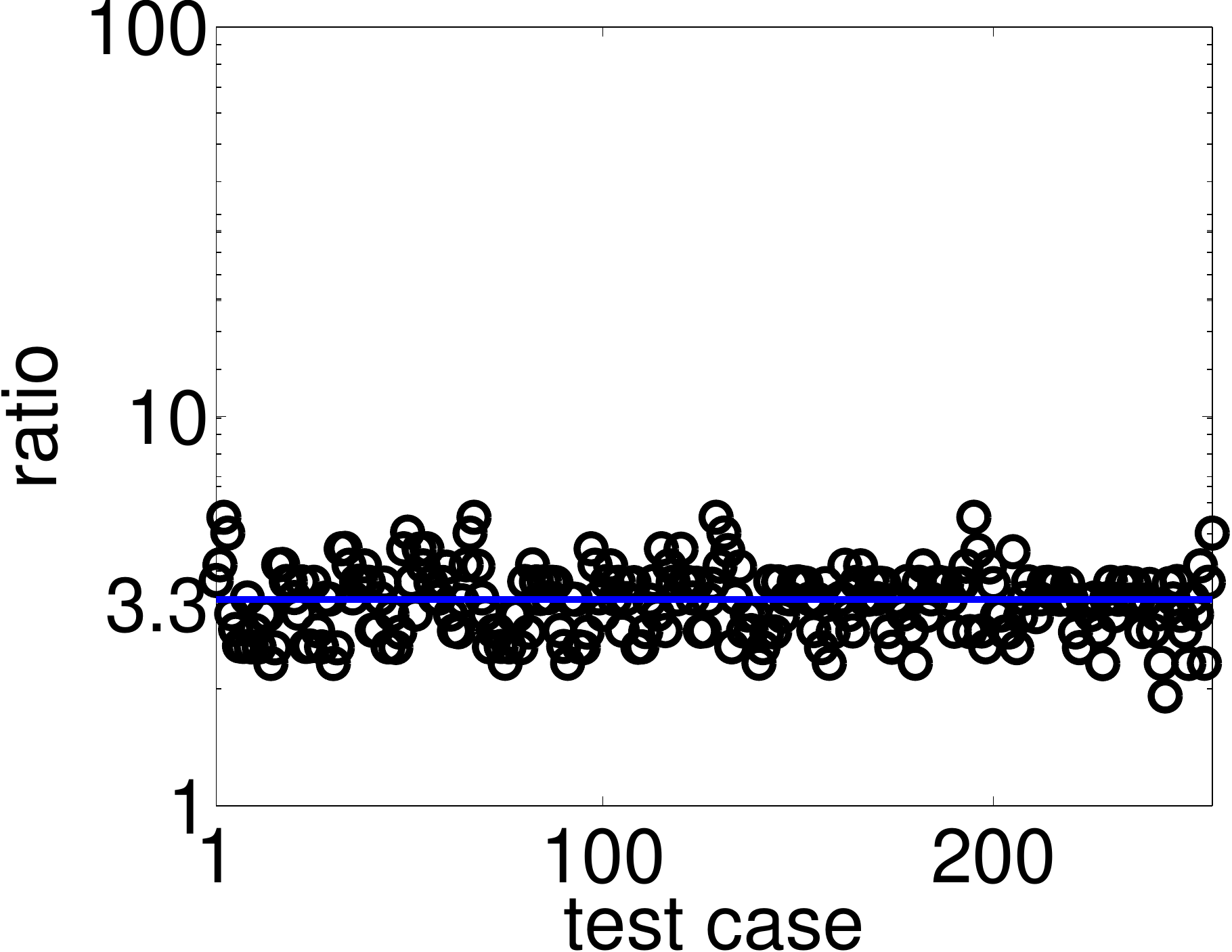} ~~
\includegraphics[width=0.4\textwidth]{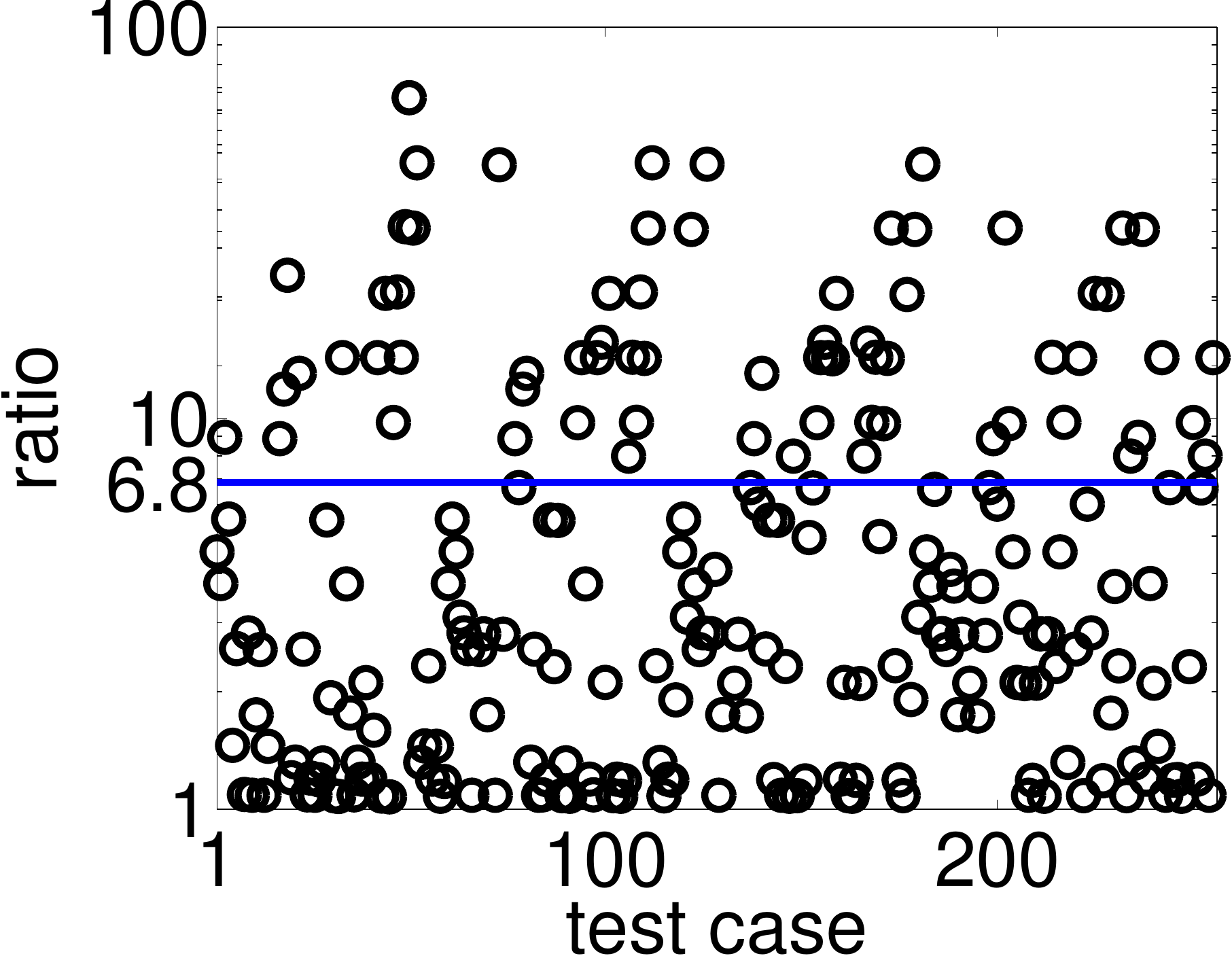}
\caption{The ratio of maximum time steps (stabilized scheme $S1$ (left) and $S2$ (right) divided by explicit scheme) for the different test cases.
The straight line shows the average ratio. The vertical axis uses log scaling.}
\label{fig:stability2}
\end{figure}

But also the two schemes employing the stabilizing terms $S_1$ and $S_2$ (Eqs. \eqref{S_st} and \eqref{S2}, resp.) perform very well. 
We divide the maximum time step of the stabilized schemes by the maximum time step of the simple explicit scheme. 
Fig. \ref{fig:stability2} shows this ratio  for the 256 test cases with $h=2\epsilon$. \\
Using $S_1$ increases the maximum time step by a factor of 1.9-5.5 (average: 3.3).
Much more diverse ratios are obtained when using $S_2$, which increases the maximum time step by a factor between 1.0 and 65.9 (average: 6.8).
Note, that the stabilizing effect of $S_1$ and $S_2$ depends on the factor $\omega$. For the ease of comparison, we set $\omega=0.2$ here, but adjusting it manually to the used mobility would allow even much higher time steps.


\subsection{Benchmark problem} \label{sec:benchmark}

We use the test setup from the two-phase flow benchmark of Hysing et al. \cite{Hysingetal_IJNMF_2009}.
It considers a single bubble rising in a liquid column for a period of $3.0$ time units. The benchmark scenario has also been studied with a diffuse interface model \cite{Aland_IJNMF_2011}.

We start with the same configuration as test case 1 in  \cite{Aland_IJNMF_2011} with $\epsilon=0.02$, $\mygamma=2\cdot 10^{-5}$, $\rho_1=1000, \rho_2=100, \nu_1=10, \nu_2=1, \sigma=24.5$.
First, let us confirm the accuracy of the $\theta$-scheme. 
To make the comparison computationally cheaper, we restrict the following studies to the time interval $[0,0.2]$.
We solve the system with different time step sizes for $\theta=1.0$ and $\theta=0.5$.
Table \ref{tab:benchmark_convergence} shows the final bubble position for all of these cases. We assume that the smallest time step together with $\theta=0.5$ gives the most accurate result and take this as a reference value to compute the errors of the other test cases.
The order of convergence (ROC) shows clearly first order convergence for $\theta=1.0$ and second order convergence for $\theta=0.5$ (Tab. \ref{tab:benchmark_convergence}).

\begin{table}
\centering
\subtable[$\theta=1.0$]{\begin{tabular}{l|rrr}
\hline\hline
 $\tau$ & position & error & ROC \\\hline
0.100&	0.51280317&	3.69E-03&		\\	
0.040&	0.51064693&	1.53E-03&	0.96\\
0.020&	0.50989193&	7.74E-04&	0.98\\
0.010&	0.50950797&	3.90E-04&	0.99\\
0.005&	0.50931372&	1.96E-04&	0.99\\
\hline\hline
\end{tabular}} 
~~~~~~~
\subtable[$\theta=0.5$]{\begin{tabular}{l|rrr}
\hline\hline
 $\tau$ & position & error & ROC \\\hline
0.100&	    0.50904344&	7.45E-05&		\\	
0.040&		0.50910657&	1.14E-05&	2.05\\
0.020&		0.50911533&	2.66E-06&	2.10\\
0.010&		0.50911762&	3.70E-07&	2.84\\
0.005&		0.50911799&	0.00E-00&		-\\
\hline\hline
\end{tabular}
}	
\caption{Bubble position, errors and rate of convergence (ROC) for different time step sizes $\tau$. The results confirm first order convergence for $\theta=1.0$ (a) and second order convergence for $\theta=0.5$ (b).}
\label{tab:benchmark_convergence}
\end{table}

So far, it did not matter whether we use 
the standard block Gauss-Seidel coupling, the fully coupled NS-CH system, or any of the introduced stabilizing terms $S_1,S_2$. All of these methods give the same computational result as long as their inner fix-point iteration converges.
However, the number of needed iterations and the time for each iteration may vary among the methods. 
Therefore, we will next analyse the performance of the four different solution methods:
\begin{itemize}
\item {\bf explicit}: standard block Gauss-Seidel coupling (see Secs. \ref{sec:block gauss seidel}-\ref{sec:simple linearization})
\item {\bf S1}: with additional stabilization term $S_1$ from Eq. \eqref{S_st}
\item {\bf S2}: with additional stabilization term $S_2$ from Eq. \eqref{S2}
\item {\bf implicit}: fully coupled NS-CH system (Eqs. \eqref{nsch_coupled}-\eqref{nsch_coupled4})
\end{itemize}
We use the same configuration as before, with $\epsilon=0.005$, $\mygamma=5\cdot 10^{-6}$.
We will use a second order Crank-Nicolson time-stepping ($\theta=0.5$) with $\tau=0.02$ which gives comparable time discretization errors as in the original benchmark paper \cite{Aland_IJNMF_2011}

\begin{table}
\centering
\begin{tabular}{l|rrrr}
\hline\hline
			& explicit & S1 & S2 & implicit \\\hline
Iterations 	& 85  & 52  & 55 & 28  \\
CPU time (s)& 116  & 71  & 89 &  95 \\
\hline\hline
\end{tabular}
\caption{Performance of different solution methods.}
\label{tab:benchmark}
\end{table}

Table \ref{tab:benchmark} shows the total number of iterations and CPU time in seconds for the different simulations.
One can see that the fully coupled system needs by far the least iterations followed by the two stabilized schemes {\bf S1},{\bf S2}.
The number of iterations for the {\bf explicit} method is more than three times higher than for the {\bf implicit} system.
However, these differences in iterations are not directly reflected in the CPU timings. 
Apparently one solution of the {\bf implicit} system is almost three times as expensive as solving the {\bf explicit} system, reducing the advantage of the {\bf implicit} system significantly. 

Also the stabilized systems {\bf S1}, {\bf S2} reduce the number of iterations compared to the {\bf explicit} scheme, while
 the resulting linear systems can be solved as fast as the {\bf explicit} system.
Consequently, the methods {\bf S1} and {\bf S2} perform best in terms of CPU time.
Method {\bf S1} is the fastest.

\begin{figure}
\includegraphics[width=0.5\textwidth]{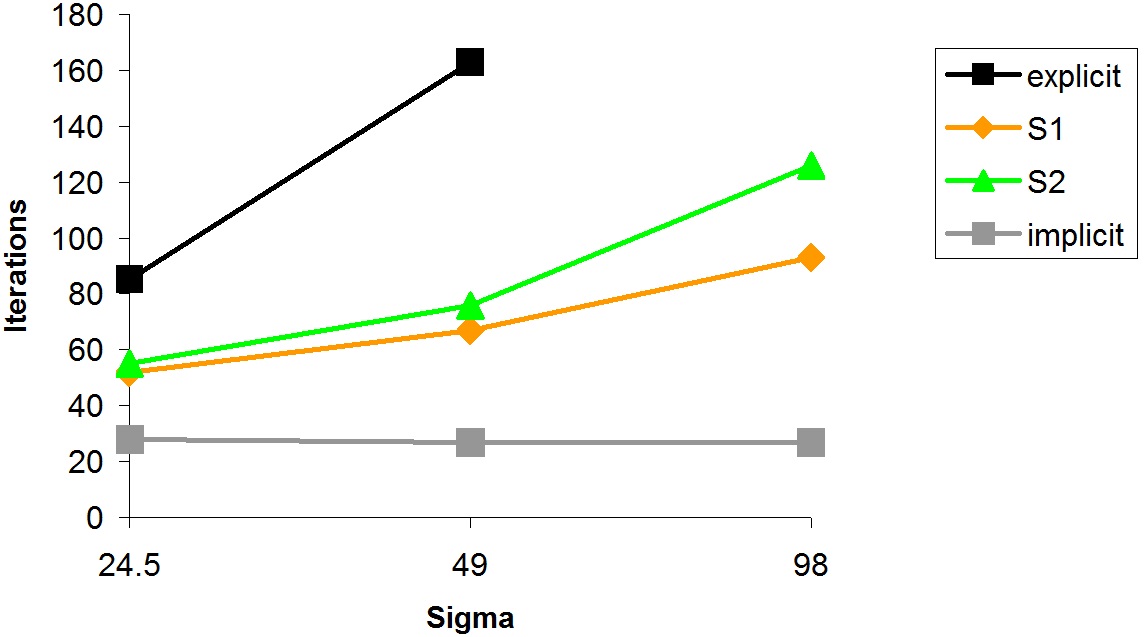}
\includegraphics[width=0.5\textwidth]{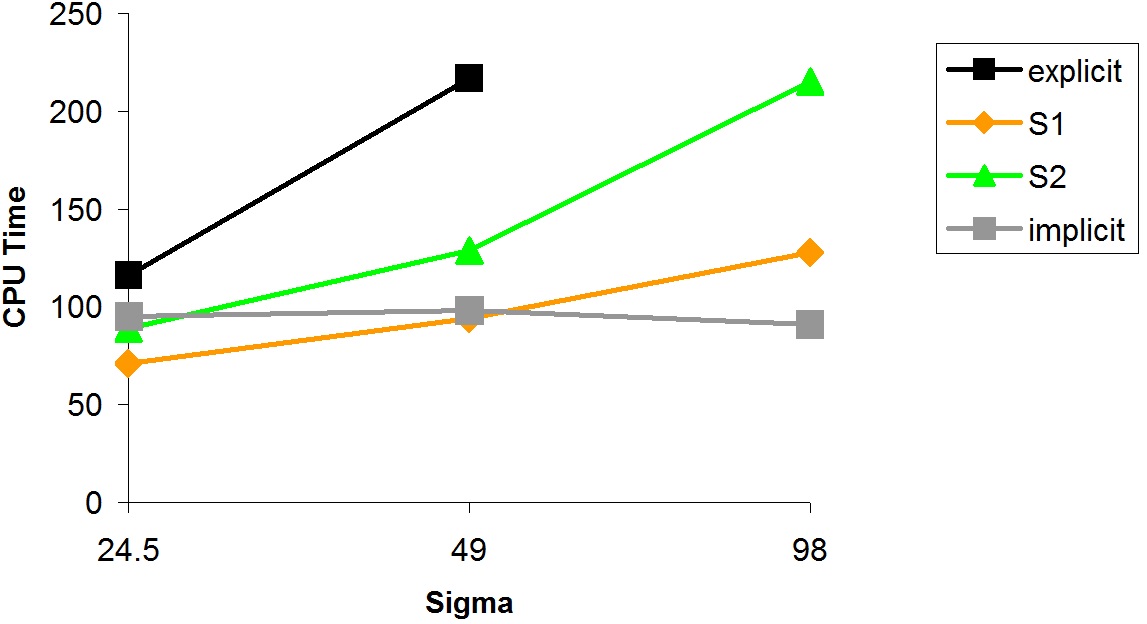}
\caption{Number of iterations and CPU times for different surface tensions $\sigma$. 
While the explicit scheme does not even converge for large surface tensions, the implicit scheme seems unaffected.}
\label{fig:benchmark1}
\end{figure}

\begin{figure}
\includegraphics[width=0.5\textwidth]{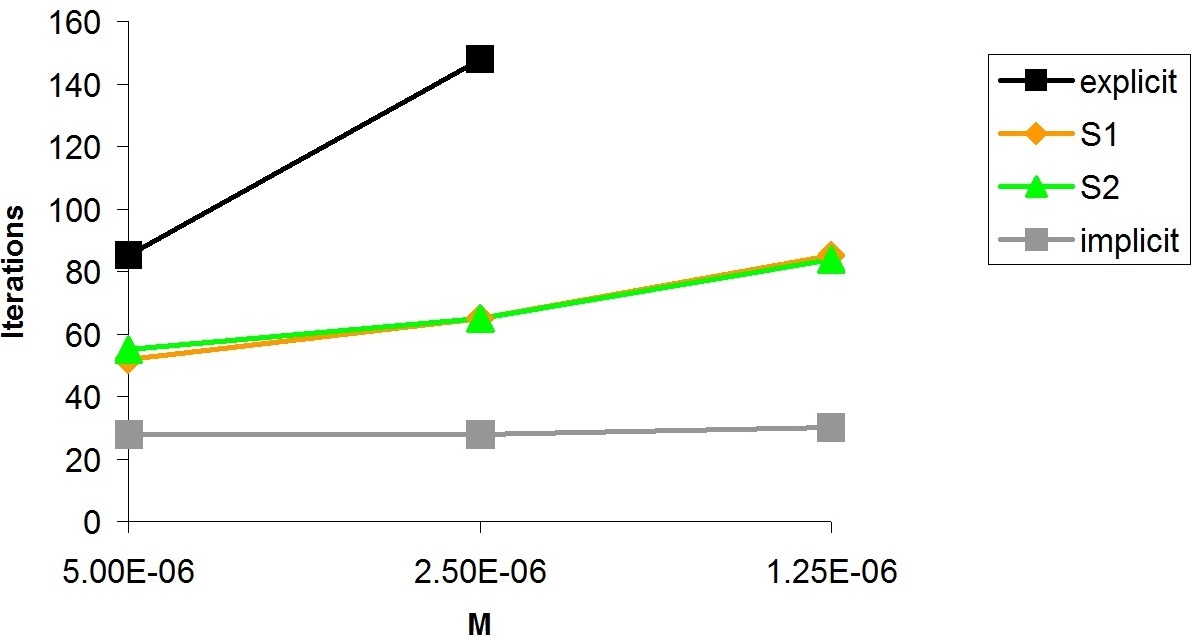}
\includegraphics[width=0.5\textwidth]{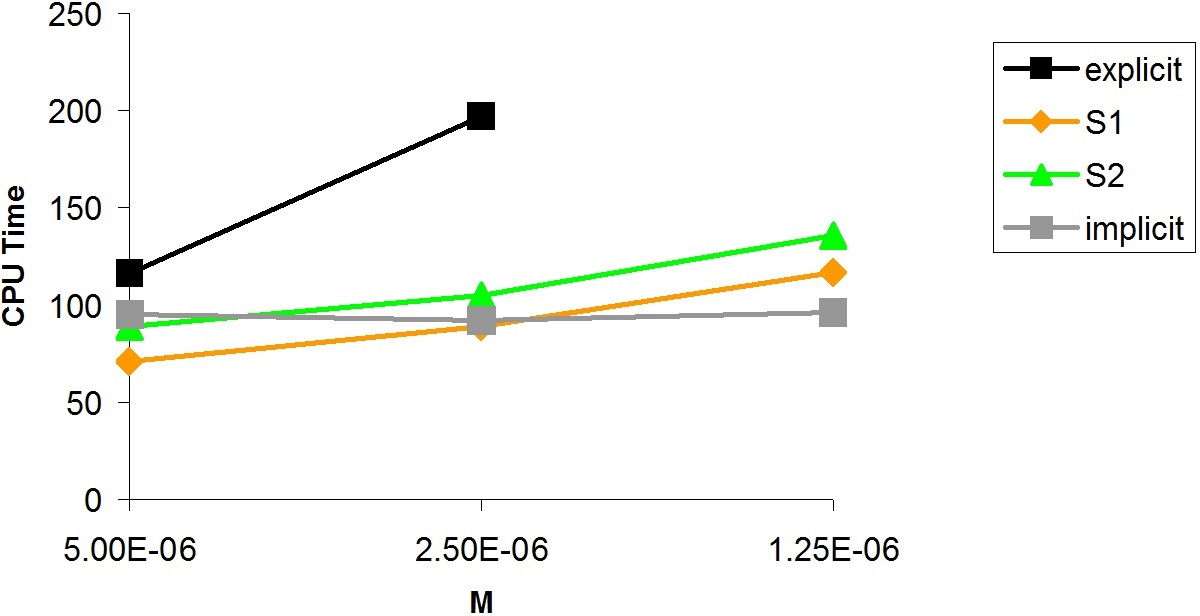}
\caption{Number of iterations and CPU times for different mobilities $M$. 
While 	the explicit scheme does not even converge for low mobility, the implicit scheme seems unaffected.}
\label{fig:benchmark2}
\end{figure}

In a next comparison we include higher surface tensions $\sigma$ and lower mobilities $\mygamma$, since we expect our improved methods to be particularly fast in these cases. Figure \ref{fig:benchmark1} shows the number of iterations and CPU times for surface tensions increased by a factor of $2$ and $4$. The number of iterations and CPU time increases for the {\bf explicit} method with increasing $\sigma$. For $\sigma=98.0$ the fix point iteration does not converge anymore.
Also the stabilized methods {\bf S1},{\bf S2} become slower for increased $\sigma$ but they still converge.
Remarkably, the {\bf implicit} method is not affected by the increase in $\sigma$, neither the number of iterations nor the CPU time changes significantly.
A very similar picture can be seen when the mobility is varied. 
In Fig.\ref{fig:benchmark2} the mobility $\mygamma$ is decreased by a factor of 2 and 4. 
Again the explicit method gets very slow and does not even converge for the smallest $\mygamma$, whereas the implicit method remains almost unaffected.


\subsection{Application to Taylor-Flow} \label{sec:taylor}
We consider a Taylor flow simulation to further demonstrate the efficiency of the proposed linearization techniques.
Taylor flow is the flow of a single elongated bubble through a narrow channel. Soon after its injection, the bubble assumes a quasi-stationary state, i.e. a fixed shape which is only advected in the direction of the channel.
Many technical applications involve Taylor flow, e.g. catalytic converters \cite{Williams01},
monolith reactors \cite{Kreutzer05}, or microfluidic
channels \cite{Muradoglu07}. In these applications, bubbles of identical
size, shape, and distance to each other are typically required. Thus, the
hydrodynamics that lead to a perfect quasi-stationary bubble are of interest in these
research areas. 
This makes Taylor flow also interesting as a general benchmark for two phase flow methods and indeed such a benchmark has been recently defined 
in 2D \cite{Alandetal_IJNMF_2013} as well as in 3D \cite{Marschall_CF_2013}.
	
\subsubsection{2D Taylor-Flow}
As in the above-mentioned benchmarks, we want to compute the quasi-stationary state of a Taylor bubble driven through the channel by a pressure difference between inlet and outlet. An implicit Euler scheme is expected to converge fast to the stationary solution due to its high numerical dissipativity.
Therefore, we use the scheme proposed  in Sec. \ref{sec:semi implicit Euler} (i.e. $\theta=1$ and only one sub-iteration) in four different versions:
\begin{itemize}
\item {\bf explicit}: no stabilization (exactly as in Sec. \ref{sec:semi implicit Euler})
\item {\bf S1}: with additional stabilization term $S_1$ from Eq. \eqref{S_st}
\item {\bf S2}: with additional stabilization term $S_2$ from Eq. \eqref{S2}
\item {\bf implicit}: fully coupled NS-CH system (Eqs. \eqref{nsch_coupled}-\eqref{nsch_coupled4} also with $\theta=1$ and only one sub-iteration)
\end{itemize}

As in \cite{Alandetal_IJNMF_2013} we use a
moving frame of reference. 
Therefore we calculate the bubble velocity ${\bf u}_b$ by ${\bf u}_b=\left(\int_\Omega (1-c){\bf u} \right)/\left(\int_\Omega 1-c\right)$ and replace the velocity of the convective terms in Eqs. (\ref{nsch_semi_implicit}) and (\ref{nsch_semi_implicit4}) by $({\bf u}-{\bf u}_b)$. 
Hence the frame of reference moves with the bubble and the quasi-stationary state becomes really-stationary.
Similar to \cite{Alandetal_IJNMF_2013} we use the parameters 
$\nu_1=\nu_2=10.0, \rho_1=\rho_2=1.0, \sigma=5000, M=3.0\cdot 10^{-6}, \epsilon=0.015, \omega=1.0$

To demonstrate the efficiency of the proposed linearization schemes, we use a simple adaptive time stepping.
Our goal is to control the CFL number which can be accomplished by choosing something like
$\tau = 0.5 h/\max(|{\bf u-u}_b|)$.
Since the usage of ${\bf u-u}_b$ would not include movement due to CH diffusion we replace this term by the phase field velocity $\partial_t c/|\nabla c|$ and end up with the following time step selection:
\begin{align}
\tau^{n} = \frac{1}{2}\tau^{n-1}\cdot\max{\frac{h|\nabla c^n|}{|c^n-c^{n-1}|}}.
\end{align}
where we choose the first time step, $\tau^{0} = 10^{-4}$.

\begin{figure}
\includegraphics[width=\textwidth]{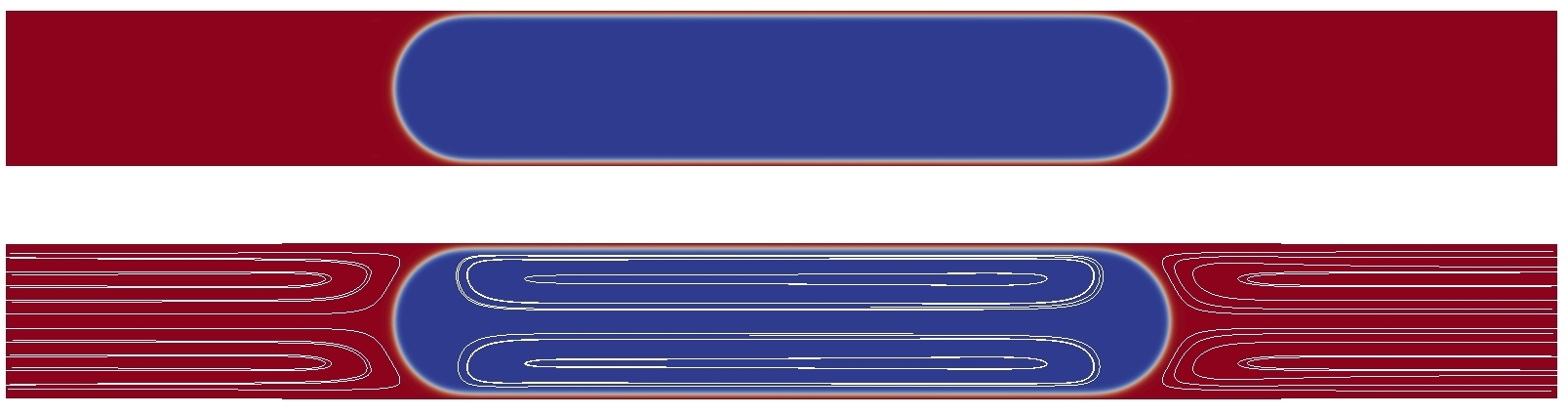}
\caption{Initial shape (top) and final shape (bottom) with flow streamlines of the Taylor bubble. The bubble shape only changed slightly during the time evolution due to high surface tension pressing the bubble against the channel wall.}
\label{fig:taylor_initial_final}
\end{figure}

The initial condition for the bubble is as in \cite{Alandetal_IJNMF_2013}:
A rod-like bubble of total length 5 is placed in the middle of a channel $\Omega=[0,10]\times[0,1]$, see Fig. \ref{fig:taylor_initial_final}.
The pressure difference between channel inflow and outflow is iteratively adjusted to give a bubble velocity ${\bf u}_b=1$.
Though it is hard to see differences between the initial and final bubble shape with the naked eye, there is some significant bubble deformation going on. 
In particular, the width of the thin liquid film between bubble and wall changes during the time evolution.

\begin{figure}
\centering
\includegraphics[width=0.9\textwidth]{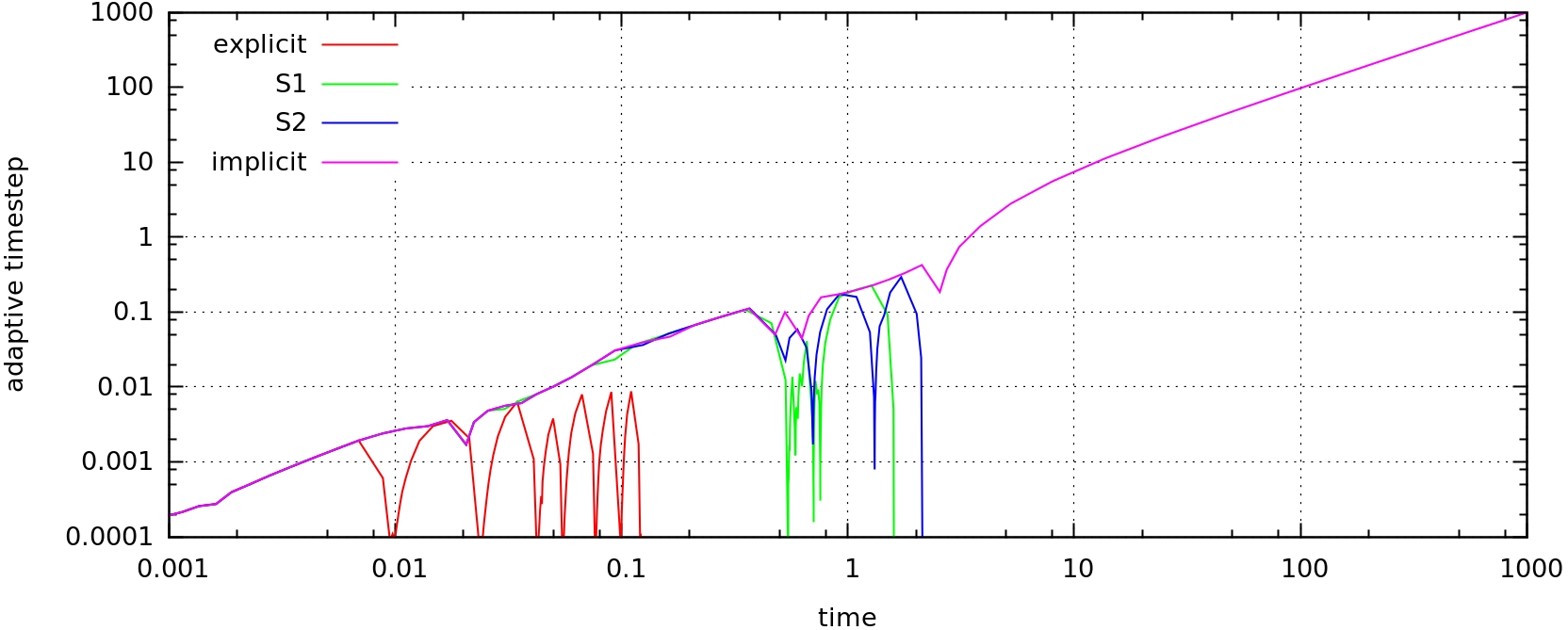} \\
\includegraphics[width=0.9\textwidth]{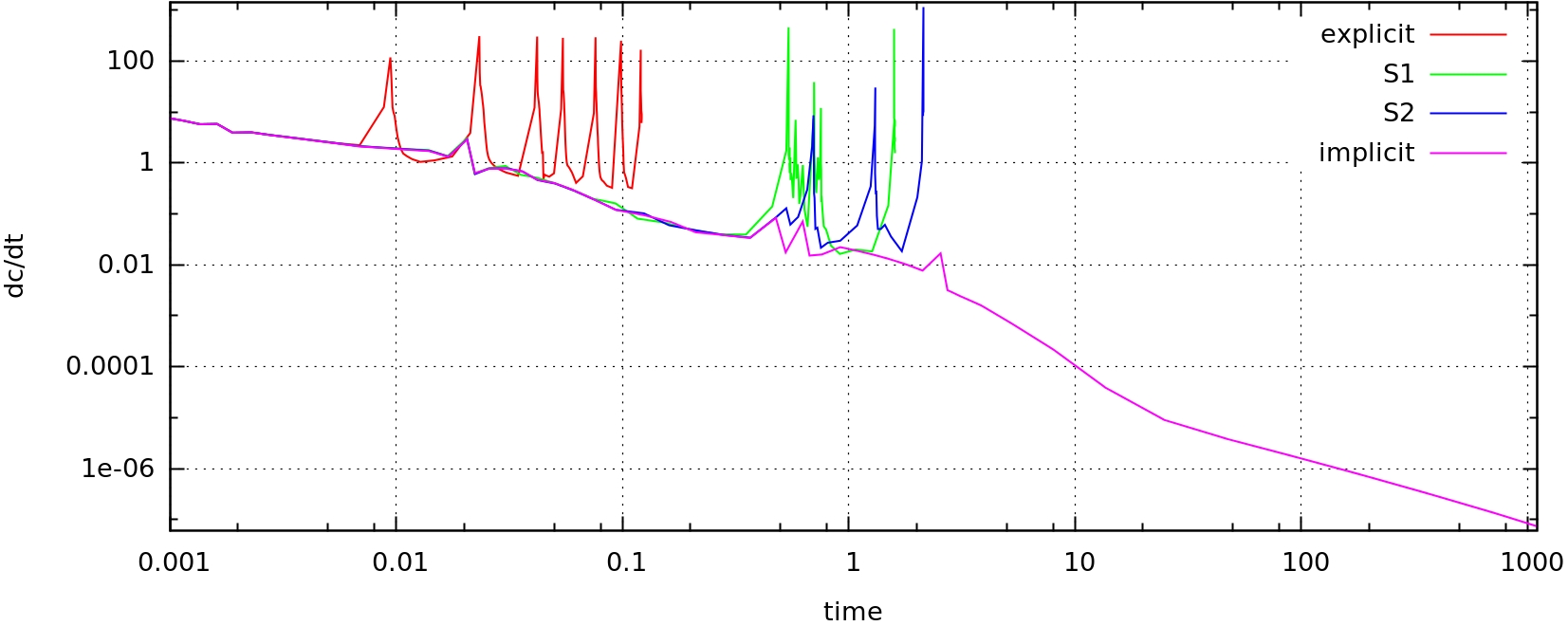}
\caption{Comparison the Taylor-flow simulation with the four methods. 
{\bf Top}: The adaptively chosen time step over time. {\bf Bottom}: The change in the phase field $|\partial_t c|$ over time. The stationary state is reached when $|\partial_t c|$ is sufficiently small.
The explicit scheme becomes unstable already for small time steps. Only the implicit scheme is stable for arbitrarily large time steps and comes very close to the stationary state.}
\label{fig:taylor_timesteps}
\end{figure}

Figure \ref{fig:taylor_timesteps} shows the adaptive time steps (top) and the change in $c$ (bottom) for the four methods. 
Let us first focus on the {\bf explicit} scheme. Although the bubble is close to its stationary shape, the {\bf explicit} scheme becomes unstable for large time steps. Once $\tau$ is larger than approximately $10^{-3}$ the interface starts to oscillate ($|\partial_t c|$ increases), which results in a rapid decrease of the time steps. After some calculations with small time steps  the bubble stabilizes and the time steps increase and the whole process starts over again. This results in an up and down of time steps which renders the explicit scheme almost useless. To reach the stationary state one would have to limit $\tau$ from above to approximately $10^{-3}$, which would require a hundred thousand of total time steps to reach a sufficiently stationary state (say at $t=100$).

A similar instability occurs for the schemes {\bf S1} and {\bf S2} but at much larger time steps, $\tau\geq 10^{-1}$.
Hence adding either of the terms $S_1$ and $S_2$ allows to increase the time steps almost by two orders of magnitude.
To reach the stationary state one would have to limit $\tau$ from above to approximately $10^{-1}$, which would require about 1000 total time steps to reach the stationary state.
Very impressing is the performance of the {\bf implicit} scheme which is stable for arbitrarily large time steps.
Thus, while the bubble approaches the stationary state, larger and larger time steps are possible. 
Hence, we obtain the stationary state in just 64 total time steps which makes absolutely worth the effort of solving NS and CH in one large system. 

We also took part in the Taylor flow benchmark paper \cite{Alandetal_IJNMF_2013} where we employed the implicit scheme. 
We could come very close to the results of sharp interface models and experimental data. 
However the sharp interface models needed several orders of magnitude more time steps to reach the same end time.
Hence, though they need a lot more degrees of freedom, diffuse interface models can be computationally cheaper and faster due to the superior possibility to couple flow equation and interface equation implicitly.

\subsubsection{3D Taylor-Flow}
Encouraged by the good results of the fully coupled time discretization for 2D-Taylor, we venture to try a 3D Taylor flow example.
We use the test setup from the 3D Taylor flow benchmark in \cite{Marschall_CF_2013}. 
Thereby a gas bubble of volume $17.5 mm^3$ is placed in a domain of size $1.98mm\times 11.88mm \times 1.98mm$. 
Exploiting the symmetry in x- and z-direction, we restrict the calculations to a quarter of this domain.
Further parameters are in the liquid phase: $\rho_L=1195.6 kg/m^3, \nu_L=28.54\cdot 10^{-3} kg/ms$, 
in the gas bubble: $\rho_G=11.95 kg/m^3, \nu_L=28.54\cdot 10^{-5} kg/ms$,
the surface tension is $\sigma= 66.69\cdot 10^{-3} kg/s^2$ and the desired final bubble velocity is $205.57 mm/s$.
For the diffuse interface model we use $M=0.9\cdot 10^{-9} m^3 s/kg$ 
 and  $\epsilon=0.03mm$. 
 
\begin{figure}
\centering
\includegraphics[width=0.38\textwidth]{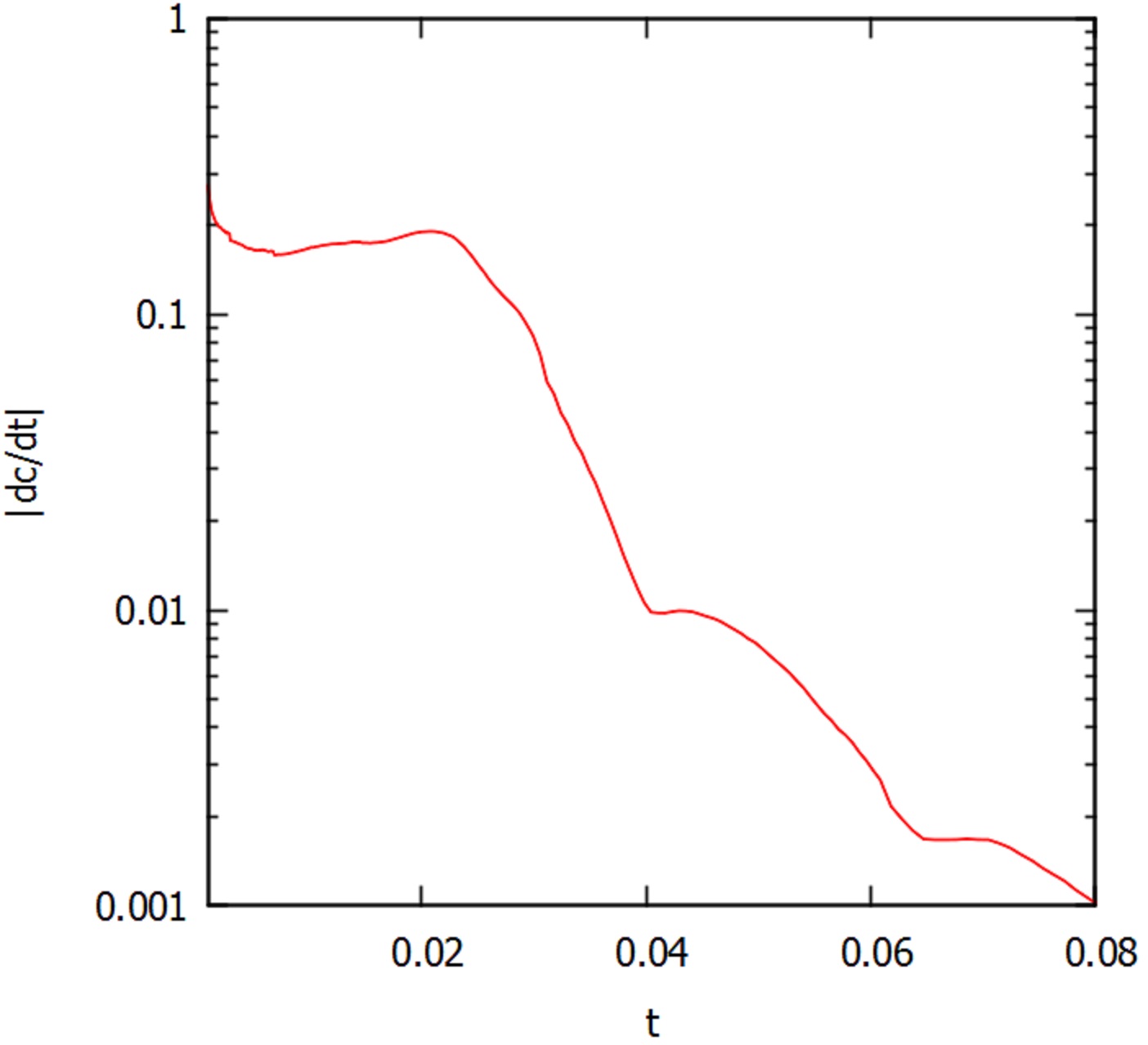} 
\includegraphics[width=0.36\textwidth]{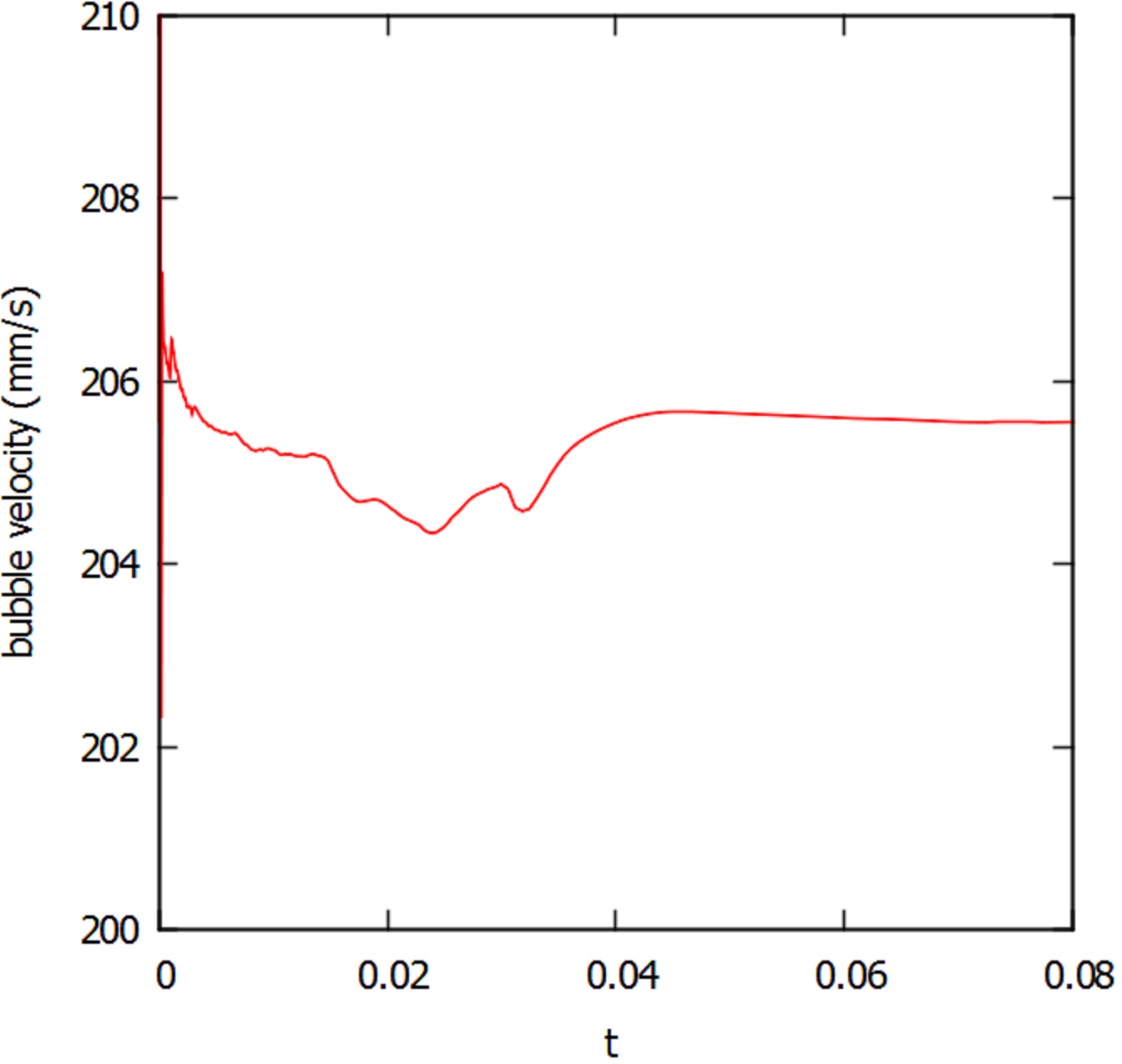} 
\caption{Time evolution of $|\partial_t c|$(left) and bubble velocity ${\bf u}_b$ (right). The change in $c$ decays exponentially (left).}
\label{fig:taylor3D_graphs}
\end{figure}

Again, we are looking for a stationary state. Hence, 
we use the fully coupled NS-CH system (Eqs. \eqref{nsch_coupled}-\eqref{nsch_coupled4} with $\theta=1$ and only one sub-iteration Starting with $\tau^0 = 2.0\cdot 10^{-5} s$, we double the time step size every 50 time steps.  
 A grid size of $h=0.087mm$ at the interface 
  gives around 1 million total degrees of freedom. 
We use an MPI based parallelization with 44 cores and the preconditioned FGMRES iteration described in Sec. \ref{sec:space discretization}.   
  One time step is solved in approximately 1 minute.
Figure \ref{fig:taylor3D_graphs} shows the time evolution of $|\partial_t c|$ and the bubble velocity ${\bf u}_b$. The change in $c$ decays exponentially.
We assume that we are sufficiently close to the stationary solution if $|\partial_t c|<10^{-3}$.

\begin{figure}
\centering
\includegraphics[width=0.95\textwidth]{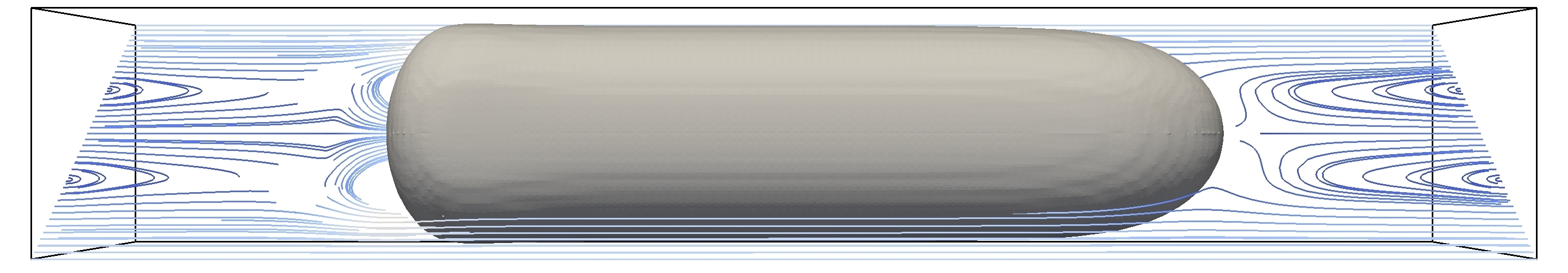} 
\caption{Quasi-stationary state of the Taylor bubble moving from the left to the right. The streamlines show a recirculating flow pattern at the front and the rear of the bubble.}
\label{fig:taylor3D}
\end{figure}

\begin{figure}
\centering
\includegraphics[width=0.48\textwidth]{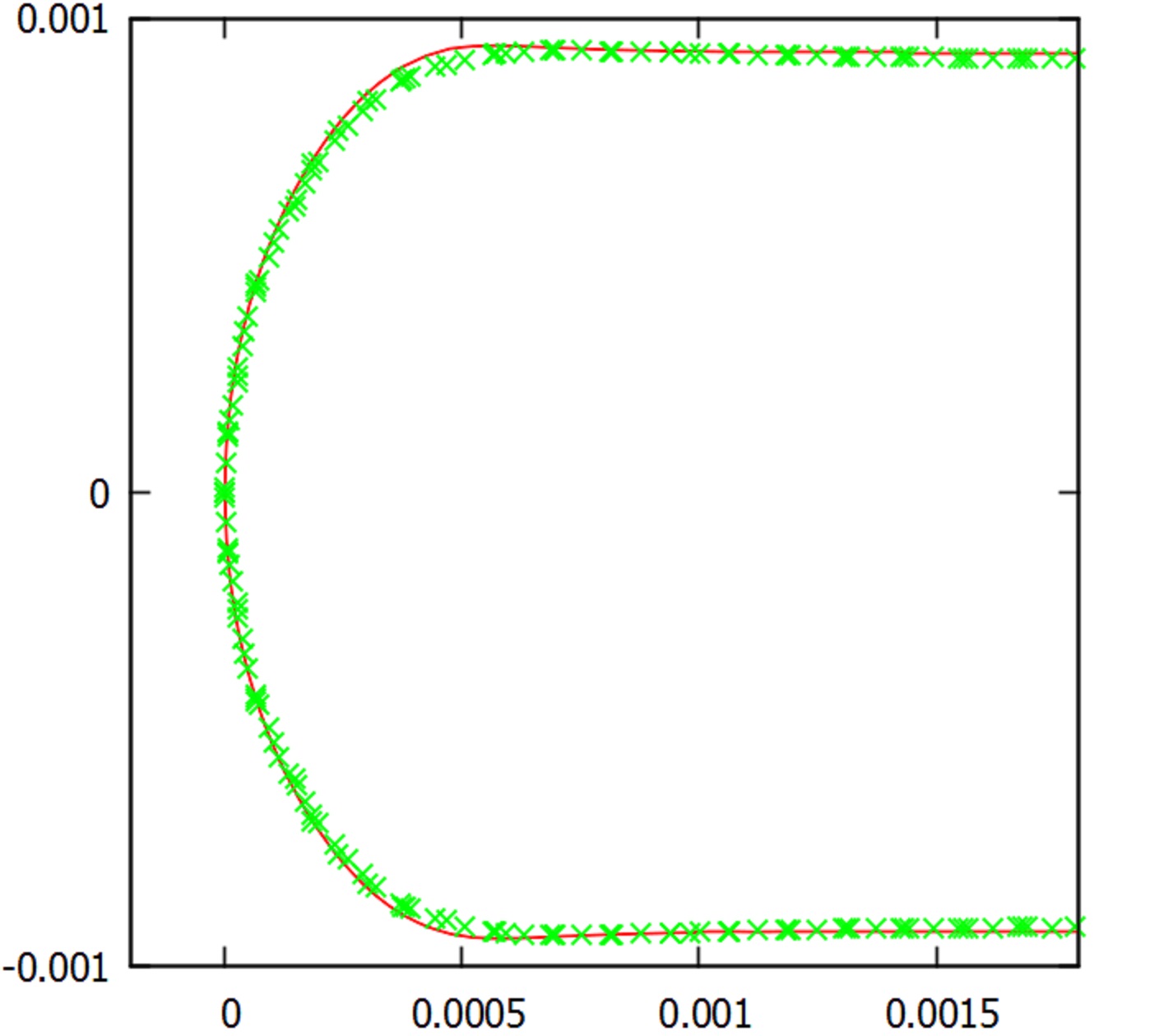} 
\includegraphics[width=0.48\textwidth]{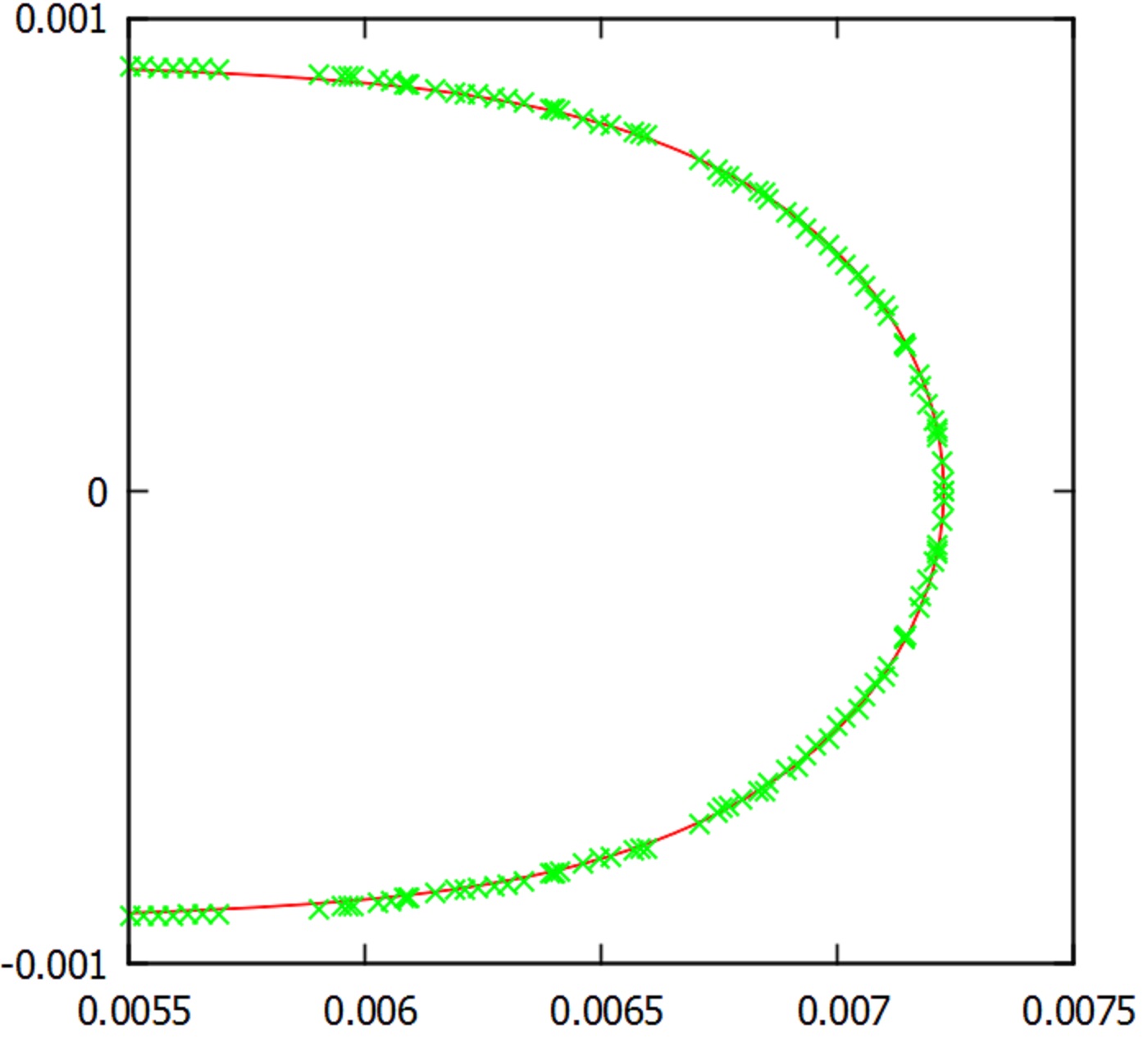} 
\caption{Final bubble shape at the rear(left) and the front(right) of the bubble. The line shows the reference solution from DROPS\cite{Marschall_CF_2013}. The crosses mark our diffuse interface solution. Both solutions are aligned at the ends.}
\label{fig:taylor3D_shape}
\end{figure}

The final bubble is depicted in Fig. \ref{fig:taylor3D} together with the streamlines showing a recirculating flow pattern at the front and the rear of the bubble. To compare the bubble shape with the reference solution from \cite{Marschall_CF_2013}, we cut the phase field vertically through the middle of the domain. The zero level set of the phase field along this slice is compared to the reference solution from DROPS \cite{Marschall_CF_2013} in Fig. \ref{fig:taylor3D_shape} and shows a quite satisfactory agreement. 
Note that these results are based on the implicit discretization scheme which here allows around ??  times larger time steps than a sequential coupling of interface and flow equation.

\section{Conclusions} \label{sec:conclusions}

In this paper, we addressed time integration strategies for the diffuse interface model for two-phase flows. 
We proposed a variant of the $\theta$-scheme together with three new linearization techniques for the surface tension.
These involve either additional stabilizing terms ($S_1$ or $S_2$), or a fully implicit coupling of the NS and CH equation.

As in all two phase flow methods the coupling between the flow and the interface equation plays a crucial role and limits the stability and the range of applicable time steps significantly.
This is particularly true if high surface tensions or small length scales are considered. 
In the common case that interface and flow equation are coupled explicitly, we could show a time step restriction of the form 
\begin{align}
\tau < C \epsilon ~\sigma^{-1/3}~ \mygamma^{1/3}~ \rho^{2/3}
\end{align}
which is very different to other two-phase flow models and in particular is independent of the grid size.
As in other two-phase flow models, this restriction can make computations extremely costly. 
Even in cases when the interface is almost stationary, too large time steps will lead to oscillations and finally destruction of the interface. 

We also showed that the proposed stabilization techniques could lift the above time step restriction.
The simple stabilizing terms $S_1$ and $S_2$ may allow an increase in time step size of about two orders of magnitude.
If a fix point or Newton sub-iteration is used in each time step, these terms can reduce the number of iterations significantly, while not affecting the accuracy. 

Very impressing is the performance of the fully implicit scheme which is stable for arbitrarily large time steps.
We demonstrate in a Taylor flow application that this superior coupling between flow and interface equation can render diffuse interface models even computationally cheaper and faster than sharp interface models.

Apart from increasing the computational performance, the improved time integration schemes allow to choose lower CH mobility, 
which may come closer to physically correct values.
Hence, the mean curvature flow included in the CH equation is suppressed which may lead to more accurate computational results.

\noindent
{\bf Acknowledgements} We acknowledge support from the German Science Foundation through grant SPP-1506 (AL 1705/1-1) and support of computing time at JSC at FZ J\"ulich.

\bibliographystyle{plain}
\bibliography{library}

\begin{thebibliography}{10}

\bibitem{Abelsetal_WS_2012}
Helmut Abels, Harald Garcke, and G{\"u}nther Gr{\"u}n.
\newblock Thermodynamically consistent, frame indifferent diffuse interface
  models for incompressible two-phase flows with different densities.
\newblock {\em Mathematical Models and Methods in Applied Sciences}, 22(03),
  2012.

\bibitem{Alandetal_IJNMF_2013}
S.~Aland, S.~Boden, A.~Hahn, F.~Klingbeil, M.~Weismann, and S.~Weller.
\newblock Quantitative comparison of {Taylor} flow simulations based on sharp-
  and diffuse-interface models.
\newblock {\em Int. J. Numer. Meth. Fluids}, 2013.

\bibitem{Aland2011_bijels}
S.~Aland, J.~Lowengrub, and A.~Voigt.
\newblock {A continuum model of colloid-stabilized interfaces}.
\newblock {\em Physics of Fluids}, 23(6):062103, 2011.

\bibitem{Aland_IJNMF_2011}
S~Aland and A~Voigt.
\newblock Benchmark computations of diffuse interface models for
  two-dimensional bubble dynamics.
\newblock {\em International Journal for Numerical Methods in Fluids},
  69(3):747--761, 2012.

\bibitem{Andersonetal_ARFM_1998}
D.M. Anderson, G.B. McFadden, and A.A. Wheeler.
\newblock {Diffuse interface methods in fluid mechanics}.
\newblock {\em Ann. Rev. Fluid Mech.}, 30:139--165, 1998.

\bibitem{Badalassi_JCP_2003}
VE~Badalassi, HD~Ceniceros, and Sanjoy Banerjee.
\newblock Computation of multiphase systems with phase field models.
\newblock {\em Journal of Computational Physics}, 190(2):371--397, 2003.

\bibitem{Boyanova_CMAM_2012}
Petia Boyanova, Minh Do-Quang, and Maya Neytcheva.
\newblock Efficient preconditioners for large scale binary cahn-hilliard
  models.
\newblock {\em Computational Methods in Applied Mathematics}, 12(1):1--22,
  2012.

\bibitem{Brackbilletal_JCP_1992}
JU~Brackbill, Douglas~B Kothe, and C1~Zemach.
\newblock A continuum method for modeling surface tension.
\newblock {\em Journal of computational physics}, 100(2):335--354, 1992.

\bibitem{umfpack}
Timothy~A. Davis.
\newblock {Algorithm 832: UMFPACK V4.3---an unsymmetric-pattern multifrontal
  method}.
\newblock {\em ACM Trans. Math. Softw.}, 30(2):196--199, June 2004.

\bibitem{Dingetal_JCP_2007}
H.~Ding, P.D.M. Spelt, and C.~Shu.
\newblock {Diffuse interface model for incompressible two-phase flows with
  large density ratios}.
\newblock {\em J. Comput. Phys.}, pages 2078--2095, 2007.

\bibitem{Do-QuangAmberg_PF_2009}
Minh Do-Quang and Gustav Amberg.
\newblock {The splash of a solid sphere impacting on a liquid surface:
  Numerical simulation of the influence of wetting}.
\newblock {\em Physics of Fluids}, 21, 2009.

\bibitem{Dziuk_NM_1990}
Gerhard Dziuk.
\newblock An algorithm for evolutionary surfaces.
\newblock {\em Numerische Mathematik}, 58(1):603--611, 1990.

\bibitem{Emmerich_AP_2008}
H.~Emmerich.
\newblock {Advances of and by phase-field modeling in condensed-matter
  physics}.
\newblock {\em Adv. Phys.}, 57:1--87, 2008.

\bibitem{Feng_SIAMJNA_2006}
X.~Feng.
\newblock {Fully Discrete Finite Element Approximations of the
  {N}avier--{S}tokes--{C}ahn-{H}illiard Diffuse Interface Model for Two-Phase
  Fluid Flows}.
\newblock {\em SIAM J. Numer. Anal.}, 44:1049--1072, 2006.

\bibitem{GruenKlingbeil_arXiv_2012}
G{\"u}nther Gr{\"u}n and Fabian Klingbeil.
\newblock Two-phase flow with mass density contrast: stable schemes for a
  thermodynamic consistent and frame-indifferent diffuse-interface model.
\newblock {\em arXiv preprint arXiv:1210.5088}, 2012.

\bibitem{Hysing_IJNMF_2006}
S~Hysing.
\newblock A new implicit surface tension implementation for interfacial flows.
\newblock {\em International Journal for Numerical Methods in Fluids},
  51(6):659--672, 2006.

\bibitem{Hysingetal_IJNMF_2009}
S.~Hysing, S.~Turek, D.~Kuzmin, N.~Parlini, E.~Burman, S.~Ganesan, and
  L.~Tobiska.
\newblock {Quantitative benchmark computations of two-dimensional bubble
  dynamics}.
\newblock {\em Int. J. Numer. Meth. Fluids}, 60:1259--1288, 2009.

\bibitem{Jacqmin_JCP_1999}
D.~Jaqmin.
\newblock {Calculation of two-phase {N}avier-{S}tokes flows using phase-field
  modelling}.
\newblock {\em J. Comput. Phys.}, 155:96--127, 1999.

\bibitem{KayWelford_SIAMJSC_2007}
D.~Kay and R.~Welford.
\newblock {Efficient Numerical Solution of {C}ahn-{H}illiard-{N}avier-{S}tokes
  Fluids in 2D}.
\newblock {\em SIAM J. Sci. Comput.}, 29:2241--2257, 2007.

\bibitem{KimLowengrub_IFB_2005}
J.~Kim and J.~Lowengrub.
\newblock Phase field modeling and simulation of three-phase flows.
\newblock {\em Interfaces and Free Boundaries}, 7:435--466, 2005.

\bibitem{Kim_JCP_2005}
Junseok Kim.
\newblock A continuous surface tension force formulation for diffuse-interface
  models.
\newblock {\em Journal of Computational Physics}, 204(2):784--804, 2005.

\bibitem{Kreutzer05}
Michiel~T. Kreutzer, Freek Kapteijn, Jacob~A. Moulijn, and Johan~J. Heiszwolf.
\newblock Multiphase monolith reactors: Chemical reaction engineering of
  segmented flow in microchannels.
\newblock {\em Chemical Engineering Science}, 60(22):5895--5916, 2005.

\bibitem{Marschall_CF_2013}
H.~Marschall, S.~Boden, C.~Lehrenfeld, C.~Falconi, U.~Hampel, A.~Reusken,
  M.~W{\"o}rner, and D.~Bothe.
\newblock {Validation of Interface Capturing and Tracking Techniques with
  different Surface Tension Treatments against a Taylor Bubble Benchmark
  Problem}.
\newblock {\em Computers and Fluids}, submitted 2013.

\bibitem{Muradoglu07}
M.~{Muradoglu}, A.~{G{\"u}nther}, and H.~A. {Stone}.
\newblock {A computational study of axial dispersion in segmented gas-liquid
  flow}.
\newblock {\em Physics of Fluids}, 19(7), 2007.

\bibitem{SingerSinger_RPP_2008}
I.~Singer-Loginova and H.~Singer.
\newblock {The phase field technique for modeling multiphase materials}.
\newblock {\em Rep. Prog. Phys.}, 71:106501, 2008.

\bibitem{Turek1999}
Stefan Turek.
\newblock {\em Efficient Solvers for Incompressible Flow Problems.: An
  Algorithmic and Computional Approach.}, volume~6.
\newblock Springer Verlag, 1999.

\bibitem{amdis}
S.~Vey and A.~Voigt.
\newblock Amdis: adaptive multidimensional simulations.
\newblock {\em Computing and Visualization in Science}, 10(1):57--67, March
  2007.

\bibitem{VillanuevaAmberg_IJMP_2006}
W.~Villanueva and G.~Amberg.
\newblock {Some generic capillary-driven flows}.
\newblock {\em Int. J. Multiphase Flow}, 32(9):1072--1086, September 2006.

\bibitem{Williams01}
J.~L. Williams.
\newblock Cheminform abstract: Monolith structures, materials, properties and
  uses.
\newblock {\em ChemInform}, 33(11), 2002.

\end{thebibliography}

\end{document}